\documentclass[a4,11pt]{article}
\usepackage[latin1]{inputenc}
\usepackage{epsfig}
\usepackage{amsmath}
\usepackage{amssymb}

\newtheorem{theo}{Theorem}[section]

\newtheorem{lema}[theo]{Lemma}
\newtheorem{coro}{Corollary}[theo]

\def\qed{\rule{1.0ex}{1.0ex} \medskip \medskip}

\def \dsp {\displaystyle}

\def\downbar#1{
\setbox10=\hbox{$#1$}
            \dimen10=\ht10 \advance\dimen10 by 2.5pt
            \ifdim \dimen10<15pt 
               \advance\dimen10 by -0.5pt
               \dimen11=\dimen10
               \advance\dimen10 by 2.5pt
               \lower \dimen11
            \else \lower \ht10 \fi
            \hbox {\hskip 1.5pt \vrule height \dimen10 depth \dp10}}
\def\upbar#1{
\setbox10=\hbox{$#1$}
            \dimen10=\ht10 \advance\dimen10 by \dp10 \advance\dimen10 by 2.5pt
            \ifdim \dimen10<15pt 
               \advance\dimen10 by 2pt \fi
            \raise 2.5pt \hbox {\hskip -1.5pt \vrule height \dimen10}}

\begin{document}

\title{ \Large{A class of orthogonal functions given by a three term recurrence formula}\thanks{This work was initiated during the exchange program CAPES(Brazil)/DGU(Spain) of 2008-2012. For this research the first and the fourth authors have also received support from CNPq and FAPESP of Brazil. The third author's research was also supported by grants from Micinn of Spain and Junta de Andaluc\'{i}a. }}
\author
{
 {C.F. Bracciali$^{a}$, J.H. McCabe$^{b}$, T.E. P\'{e}rez$^{c}$ and A. Sri Ranga$^{a}$} \\[1ex]
  {\small $^{a}$Departamento de Matem\'{a}tica Aplicada, IBILCE, }\\
  {\small  UNESP - Universidade Estadual Paulista} \\
  {\small 15054-000, São José do Rio Preto, SP, Brazil }\\[1ex]
  {\small $^{b}$Department of Applied Mathematics, School of Mathematics,} \\
  {\small  University of St.Andrews, Scotland }\\[1ex]
    {\small $^{c}$Departamento de Matem\'{a}tica Aplicada, } \\
  {\small  Universidad de Granada, 18071 Granada, Spain }
}


\maketitle

\thispagestyle{empty}

\begin{abstract}
   The main goal in this manuscript is to present a class of functions satisfying a certain orthogonality property for which there also exists a three term recurrence formula. This class of functions, which can be considered as an extension to the class of symmetric orthogonal polynomials on $[-1,1]$, has a complete connection to the  orthogonal polynomials on the unit circle.  Quadrature rules and other properties based on the zeros of these functions are also considered.
\end{abstract}


\setcounter{equation}{0}
\section{Introduction  }

 Let  $\Omega_{m}$ be the linear space of ``real'' functions on $[-1,1]$ defined as follows. \\[-1ex]

  $\Omega_0 \equiv \mathbb{P}_0$ \ and \ $\Omega_m$  for $m \geq 1$ is such that  if \ $\mathcal{F} \in \Omega_{m}$ \ then  \ $\mathcal{F}(x) = B^{(0)}(x) + \sqrt{1-x^2}\, B^{(1)}(x)$, where $B^{(0)}(x) \in \mathbb{P}_m$ and $B^{(1)}(x) \in \mathbb{P}_{m-1}$ satisfy
\[
  \begin{array}l
    B^{(0)}(-x) = (-1)^{m}B^{(0)}(x) \ \ \ \mbox{and}  \ \ \ B^{(1)}(-x) = (-1)^{m-1}B^{(1)}(x).
  \end{array}
\]
Here, $\mathbb{P}_m$ represents the linear space of real polynomials of degree at most $m$.

This means, if $\mathcal{F} \in \Omega_{2n}$ then $B^{(0)}$ is an even polynomial of degree at most $2n$ and $B^{(1)}$ is an odd polynomial of degree at most $2n-1$.  Likewise, if $\mathcal{F} \in \Omega_{2n+1}$ then $B^{(0)}$ is an odd polynomial of degree at most $2n+1$ and $B^{(1)}$ is an even polynomial of degree at most $2n$.

Note that the dimension of $\Omega_m$ is $m+1$.  As  an example of a basis for $\Omega_{2n}$ we have
\[
   \{ 1, \ x\sqrt{1-x^2}, \ x^2, \   \ldots, \ x^{2n-1}\sqrt{1-x^2}, \ x^{2n} \}
\]
and as an example of a basis for $\Omega_{2n+1}$ we can state
\[
   \{ \sqrt{1-x^2}, \ x, \ x^2\sqrt{1-x^2}, \  \ldots, \ x^{2n}\sqrt{1-x^2}, \ x^{2n+1}\}.
\]

Our interest in the studies of such functions are many.  Apart from the interesting properties such as three term recurrence formula, orthogonality and  quadrature formulas that can be associated with these functions as shown in this manuscript, solutions of the following differential equations, with integer $m$, are also of these type of functions (see \cite{DimRan-2013}).
\[ (1-x^2)\mathcal{F}^{\prime\prime}(x) - \big[(2\lambda+1)x - 2 \eta \sqrt{1-x^2}\,\big] \mathcal{F}^{\prime}(x) +
      m\big[m + 2 \lambda +  \frac{2 m \eta x}{\sqrt{1-x^2}}\big] \mathcal{F}(x) = 0.
\]
When $\eta = 0$, the solutions of the above differential equations are the ultraspherical polynomials (see, for example,  \cite{Chihara-book, Szego-book}).

We will also describe in section \ref{Sec-OPUC-Connection} of this manuscript the connection that the functions considered here have with orthogonal polynomials on the unit circle (OPUC).

For $\mathcal{F}(x) = B^{(0)}(x) + \sqrt{1-x^2}\, B^{(1)}(x) \in \Omega_m$, by setting
\[
  B^{(0)}(x) = \sum_{j=0}^{\lfloor m/2 \rfloor} b_{2j}^{(0)}\,x^{m-2j} \ \ \mbox{and} \ \ B^{(1)}(x) = \sum_{j=0}^{\lfloor (m-1)/2 \rfloor}  b_{2j}^{(1)}\,x^{m-1-2j},
\]
we say that the function $\mathcal{F}$ is of exact degree $m$ if  $(b_{0}^{(0)})^2 + (b_{0}^{(1)})^2 > 0$. The nonnegative number $(b_{0}^{(0)})^2 + (b_{0}^{(1)})^2$ may be called the lead factor of $\mathcal{F}$. The coefficients $b_{0}^{(0)}$ and $b_{0}^{(1)}$ may be  referred to as the first and second leading coefficients of $\mathcal{F}$, respectively.

The specific aim of this manuscript is to consider some properties, in particular, the orthogonal properties, of the  sequence of functions $\{\mathcal{W}_m(x)\}$, where $\mathcal{W}_m(x) \in \Omega_m$,  given by
\begin{equation} \label{TTRR-Wm}
  \begin{array}{l}
  \mathcal{W}_0(x) = \gamma_0,\quad \mathcal{W}_1(x) = (\gamma_1 x - \beta_{1} \sqrt{1-x^2})\gamma_0,  \\[1ex]
  \mathcal{W}_{m+1}(x) = \big[\gamma_{m+1}x - \beta_{m+1} \sqrt{1-x^2}\big] \mathcal{W}_{m}(x) - \alpha_{m+1} \mathcal{W}_{m-1}(x), \quad m \geq 1.
  \end{array}
\end{equation}
Here,  $\{\alpha_{m}\}_{m=2}^{\infty}$, $\{\beta_m\}_{m=1}^{\infty}$ and  $\{\gamma_m\}_{m=0}^{\infty}$ are sequences of real numbers.

\setcounter{equation}{0}
\section{Functions in $\Omega_{m}$ and self-inversive polynomials}

It is known that a polynomial $Q_m$ of degree at most $m$ is self-inversive of degree $m$ if $z^m \overline{Q_m(1/\overline{z})} = c_m Q_m(z)$, where $|c_m| =1$. As one of the earliest references to self-inversive polynomials we cite Bonsall and Marden \cite{BonsallMarden-1952}. The most interesting self-inversive polynomials are those polynomials with all their zeros on the unit circle.  Characterizing self-inversive polynomials with  zeros on the unit circle has been of considerable interest (see, for example, \cite{JoNjTh-1989, Schinzel-2005, Sinclair-Vaaler-2008}).

In this manuscript we adopt the definition that $Q_m$ is a self-inversive polynomial of degree $m$ if $Q_m$ is a polynomial of degree at most $m$ and satisfy
\[
    Q_m^{\ast}(z) = z^{m}\overline{Q_m(1/\overline{z})} = Q_m(z).
\]
Note that we have assumed  $c_m =1$ from the original definition of self-inversive polynomials. We remark that  what is considered in the present manuscript as  self-inversive polynomials are, as in \cite{Sinclair-Vaaler-2008}, are known as conjugate reciprocal polynomials.

Functions belonging to $\Omega_m$ are connected to self-inversive polynomials of degree $m$. That is, given $\mathcal{F}_m \in \Omega_m$ then associated with it there exists a unique $Q_m$ which is a self-inversive polynomials  of degree  $m$. Precisely,  $ e^{-im\theta/2}Q_m(e^{i\theta}) = \mathcal{F}_m(x)$,  where $x = \cos(\theta/2)$.

The following lemma, which will be of considerable use in this manuscript, gives a  more precise statement regarding this connection.

\begin{lema} \label{Lemma1} Let $x=\cos(\theta/2)$. Then the polynomial $Q_m$ is  self-inversive of degree $m$ if and only if
\[
      e^{-im\theta/2} Q_m(e^{i\theta}) = \mathcal{F}_m(x) = B_m^{(0)}(x) + \sqrt{1-x^2}\, B_{m}^{(1)}(x),
\]
where $B_m^{(0)}$ and $B_{m}^{(1)}$ are real polynomials of degree at most $m$ and $m-1$, respectively, and  satisfying the symmetry
\[
     B_{m}^{(0)}(-x) = (-1)^m B_{m}^{(0)}(x) \quad \mbox{and} \quad B_{m}^{(1)}(-x) = (-1)^{m-1} B_{m}^{(1)}(x).
\]
Thus,
\[
    |Q_m(e^{i\theta})|^2  = \big[B_{m}^{(0)}(x) + \sqrt{1-x^2}\, B_{m}^{(1)}(x)\big]^2 \quad \mbox{and} \quad Q_m(1) = B_{m}^{(0)}(1).
\]
Moreover, $Q_m$ is a self-inversive polynomial with real coefficients if and only if $B_{m}^{(1)}$ is identically zero.
\end{lema}

\noindent{\bf Proof}. This Lemma has also been stated and proved  in \cite{DimIsmailRan-2012}. However, for completeness and for a better understanding of the use of this lemma, we give here a sketch of its proof.

Given any polynomial $Q_m$ of degree at most $m$, not necessarily self-inversive, by setting  $Q_m(z) = \sum_{j=0}^{m} (c_j^{(m)} + i\, d_j^{(m)}) z^j$ we can write
\[
    \begin{array}{rl}
    z^{-m+1/2} Q_{2m-1}(z) & = \dsp  \sum_{j=0}^{m-1} (c_{m-1-j}^{(2m-1)} + c_{m+j}^{(2m-1)})\frac{z^{j+1/2}+z^{-j-1/2}}{2} \   \\[3ex]
     & \dsp \qquad - \ \sum_{j=0}^{m-1}\ i(d_{m-1-j}^{(2m-1)} - d_{m+j}^{(2m-1)})\frac{z^{j+1/2}-z^{-j-1/2}}{2} \\[3ex]
     & \dsp \qquad \qquad  + \ \sum_{j=0}^{m-1} i(d_{m-1-j}^{(2m-1)} + d_{m+j}^{(2m-1)})\frac{z^{j+1/2}+z^{-j-1/2}}{2} \ \ \\[3ex]
     & \dsp \qquad \qquad \qquad - \ \sum_{j=0}^{m-1} (c_{m-1-j}^{(2m-1)} - c_{m+j}^{(2m-1)})\frac{z^{j+1/2}-z^{-j-1/2}}{2}
  \end{array}
\]
and
\[
  \begin{array}{rl}
    z^{-m} Q_{2m}(z) & = \dsp c_m^{(2m)} + \sum_{j=1}^{m} \left[(c_{m-j}^{(2m)} + c_{m+j}^{(2m)})\frac{z^{j}+z^{-j}}{2} \ - \ i(d_{m-j}^{(2m)} - d_{m+j}^{(2m)})\frac{z^{j}-z^{-j}}{2} \right] \\[3ex]
      & \dsp \qquad + \ i \,d_m^{(2m)} + \sum_{j=1}^{m} \left[i(d_{m-j}^{(2m)} + d_{m+j}^{(2m)})\frac{z^{j}+z^{-j}}{2} \ - \ (c_{m-j}^{(2m)} - c_{m+j}^{(2m)})\frac{z^{j}-z^{-j}}{2} \right] .
    \end{array}
\]
Then for $z = e^{i\theta}$, with $x  = (z^{1/2} + z^{-1/2})/2 = \cos(\theta/2)$, we obtain
\begin{equation*} \label{General-Mapping}
\begin{array}l
     e^{-im\theta/2} Q_{m}(e^{i\theta}) = \mathcal{F}_{m}(x) + i\, \tilde{\mathcal{F}}_{m}(x) \quad \mbox{and} \quad |Q_{m}(e^{i\theta})|^2 = [\mathcal{F}_{m}(x)]^2 + [\tilde{\mathcal{F}}_{m}(x)]^2,
\end{array}
\end{equation*}
where the functions $\mathcal{F}_m$ and $\tilde{\mathcal{F}}_m$, defined for $x \in [-1,1]$, satisfy
\[
\begin{array}l
    \dsp \mathcal{F}_{2m-1}(x) =  \sum_{j=0}^{m-1} (c_{m-1-j}^{(2m-1)} + c_{m+j}^{(2m-1)})T_{2j+1}(x) + (d_{m-1-j}^{(2m-1)} - d_{m+j}^{(2m-1)})\sqrt{1-x^2}\,U_{2j}(x), \\[3ex]
    \dsp \tilde{\mathcal{F}}_{2m-1}(x) =  \sum_{j=0}^{m-1} (d_{m-1-j}^{(2m-1)} + d_{m+j}^{(2m-1)})T_{2j+1}(x) - (c_{m-1-j}^{(2m-1)} - c_{m+j}^{(2m-1)})\sqrt{1-x^2}\,U_{2j}(x), \\[3ex]
    \dsp \mathcal{F}_{2m}(x) = c_{m}^{(2m)} + \sum_{j=1}^{m} (c_{m-j}^{(2m)} + c_{m+j}^{(2m)})T_{2j}(x) + (d_{m-j}^{(2m)} - d_{m+j}^{(2m)})\sqrt{1-x^2}\,U_{2j-1}(x), \\[3ex]
    \dsp \tilde{\mathcal{F}}_{2m}(x) = d_{m}^{(2m)} + \sum_{j=1}^{m} (d_{m-j}^{(2m)} + d_{m+j}^{(2m)})T_{2j}(x) - (c_{m-j}^{(2m)} - c_{m+j}^{(2m)})\sqrt{1-x^2}\,U_{2j-1}(x). \\[1ex]
\end{array}
\]
Here,
\[
T_j(x) = \cos(j\theta/2) = \frac{1}{2}(z^{j/2} + z^{-j/2}) \quad \mbox{and} \quad U_j(x) = \frac{\sin((j+1)\theta/2)}{\sin(\theta/2)} = \frac{(z^{(j+1)/2} - z^{-(j+1)/2})} {(z^{1/2} - z^{-1/2})},
\]
are respectively the Chebyshev polynomials of the first and second kind.

These relations enable one to establish a connection  between any polynomial $Q_m(z)$ and two functions $\mathcal{F}_m(x)$ and $\tilde{\mathcal{F}}_m(x)$ both in $\Omega_m$.

However, if $Q_m(z) = \sum_{j=0}^{m} (c_j^{(m)} + i\, d_j^{(m)}) z^j$ is self-inversive then $c_j^{(m)} = c_{m-j}^{(m)}$ and $d_j^{(m)} = -d_{m-j}^{(m)}$ and hence $\tilde{\mathcal{F}}_m(x)$  is identically zero.  \hfill \qed

An immediate consequence of this Lemma is the following.

\begin{theo} \label{Thm-Max-Zeros}
If $\mathcal{F} \in \Omega_m$ then the number of zeros of $\mathcal{F}$ in $[-1, 1]$ can not exceed $m$.
\end{theo}

Another interesting result regarding functions in $\Omega_m$ is their interpolation property, which is virtually the interpolation property of self-inversive polynomials on the unit circle $|z| =1$.  Results regarding interpolation by polynomials, including  the idea behind the proof of the theorem below, are well known and can be found in any numerical analysis texts. For a recent reference to such a text we cite \cite{Phillips-book-2003}.

\begin{theo} \label{Thm-Interpolation}
   Given the $m+1$ pairs of real numbers $(x_j,y_j)$, $j=1,2, \ldots, m+1$, where $-1 < x_1 < x_2 < \ldots < x_{m+1} < 1$, then there exists a unique $\mathcal{F} \in \Omega_m$ such that
\[
     \mathcal{F}(x_j) = y_j, \quad j=1,2, \ldots, m+1.
\]
Moreover, this  interpolation function can be given by
\begin{equation} \label{Eq-Interpolation}
     \mathcal{F}(x) = \sum_{k=1}^{m+1} \mathcal{L}_{k}(x)\, y_k,
\end{equation}
where
\[
   \mathcal{L}_{k}(x) = z^{-m/2} \ z_k^{m/2}\!\!
   \prod_{ \begin{array}c
            l=1 \\[-0.5ex]
            l \neq k
           \end{array}
          }^{m+1}
   \frac{z-z_l}{z_k-z_l}, \qquad k=1,2, \ldots, m+1,
\]
with $z = e^{i\theta}$, $\theta = 2 \arccos(x)$, $z_k = e^{i\theta_k}$ and  $\theta_k = 2 \arccos(x_k)$.

\end{theo}

\noindent{\bf Proof}. Uniqueness follows from Theorem \ref{Thm-Max-Zeros}. That is if $\mathcal{F}$ and $\tilde{\mathcal{F}}$ are two different functions in $\Omega_m$ such that
\[
     \mathcal{F}(x_j) = y_j \quad \mbox{and} \quad \tilde{\mathcal{F}}(x_j) = y_j, \quad j=1,2, \ldots, m+1,
\]
then $\mathcal{G}(x)$, where $\mathcal{G}(x) = \mathcal{F}(x) - \tilde{\mathcal{F}}(x) \in \Omega_m$ and $\mathcal{G}_m(x) \neq 0$, has $m+1$ zeros in $(-1,1)$ contradicting Theorem \ref{Thm-Max-Zeros}.

To show the existence we construct the required function  as follows.
It is easily seen that the scaled Lagrange polynomials
\[
     \ z_k^{m/2}\!\!
   \prod_{\begin{array}c
           l=1 \\[-0.5ex]
           l \neq k
          \end{array}}^{m+1}
   \frac{z-z_l}{z_k-z_l}, \qquad k=1,2, \ldots, m+1,
\]
defined on the set $\{z_1, z_2, \ldots, z_m\}$ are self-inversive. Consequently,
$\mathcal{L}_{k}(x)$ are in $\Omega_m$  and that they  satisfy
\[
   \mathcal{L}_{k}(x_j) = \delta_{jk}, \quad j=1,2, \ldots, m+1.
\]
Hence, the formula (\ref{Eq-Interpolation})  immediately leads to the required interpolation function. \hfill \qed

\setcounter{equation}{0}
\section{Some basic properties }

\begin{theo} \label{Thm-Basis}
Let the sequence functions $\{\mathcal{F}_m\}$, where $\mathcal{F}_0(x) = b_{0,0}^{(0)}$,
\[
  \mathcal{F}_m(x) = \sum_{j=0}^{\lfloor m/2 \rfloor}  b_{m,2j}^{(0)}\,x^{m-2j} \  + \ \sqrt{1-x^2}\, \sum_{j=0}^{\lfloor (m-1)/2 \rfloor}  b_{m,2j}^{(1)}\,x^{m-1-2j}, \quad m \geq 1,
\]
be such that $b_{m+1,0}^{(0)}b_{m,0}^{(0)}+b_{m+1,0}^{(1)}b_{m,0}^{(1)} \neq 0$, $m \geq 0$.  Here, $b_{0,0}^{(1)} = 0$. Then the following hold. \\[-0.5ex]

\noindent 1.\ \ A basis for $\Omega_{2m}$ is
\[
   \{\mathcal{F}_{2m}(x), \sqrt{1-x^2}\mathcal{F}_{2m-1}(x), \mathcal{F}_{2m-2}(x), \ldots, \mathcal{F}_2(x), \sqrt{1-x^2}\mathcal{F}_{1}(x), \mathcal{F}_{0}(x)\};
\]

\noindent 2.\ \ A basis for $\Omega_{2m+1}$ is
\[
   \{\mathcal{F}_{2m+1}(x), \sqrt{1-x^2}\mathcal{F}_{2m}(x), \mathcal{F}_{2m-1}(x), \ldots, \mathcal{F}_3(x), \sqrt{1-x^2}\mathcal{F}_{2}(x), \mathcal{F}_{1}(x), \sqrt{1-x^2}\mathcal{F}_{0}(x)\}.
\]
\end{theo}

\noindent{\bf Proof}. Since $b_{m+1,0}^{(0)}b_{m,0}^{(0)}+b_{m+1,0}^{(1)}b_{m,0}^{(1)} \neq 0$ for $m \geq 0$ implies $(b_{m,0}^{(0)})^2+(b_{m,0}^{(1)})^2 > 0$ for $m \geq 0$,  the function $\mathcal{F}_m$ is of exact degree $m$ for $m \geq 0$.

Now to prove the theorem, by observing that the dimension of $\Omega_{m}$ is $m+1$,  all we have to do is to verify that the above sets are linearly independent. We prove this for the even indices and the proof for the odd indices is similar.

Clearly $\mathcal{F}_0(x) = b_{0,0}^{(0)} \neq 0$ is a basis for $\Omega_{0}$.  Now we verify that
\[
   \{\mathcal{F}_{2}(x), \sqrt{1-x^2}\mathcal{F}_{1}(x), \mathcal{F}_{0}(x)\}
\]
is a basis for $\Omega_{2}$.

Let $c_0, c_1, c_2$ be such that  $c_0\, \mathcal{F}_{2}(x) + c_1\, \sqrt{1-x^2}\,\mathcal{F}_{1}(x) + c_2\, \mathcal{F}_{0}(x) = 0$.  Since the dimension of $\Omega_{2}$ is 3, we need to verify that this is possible only if $c_0 = c_1 = c_2=0$. By considering the coefficients of $x^2$ and $x \sqrt{1-x^2}$, we have
\[
 \begin{array}{rl}
   c_0\,b_{2,0}^{(0)} - c_1\,b_{1,0}^{(1)} &= 0, \\[1ex]
   c_0\,b_{2,0}^{(1)} + c_1\,b_{1,0}^{(0)} &= 0.
 \end{array}
\]
The determinant of this system is $b_{2,0}^{(0)}b_{1,0}^{(0)}+b_{2,0}^{(1)}b_{1,0}^{(1)}$, which is different from zero.  We must therefore have $c_0 = c_1 = 0$. This reduces our verification to finding  $c_2$ such that $c_2\,\mathcal{F}_0(x) = 0$.  Clearly $c_2=0$, which follows from $\mathcal{F}_0(x) = b_{0,0}^{(0)} \neq 0$.

Now assuming %
\begin{equation}  \label{Assumption-1}
   \{\mathcal{F}_{2m}(x), \sqrt{1-x^2}\mathcal{F}_{2m-1}(x), \mathcal{F}_{2m-2}(x), \ldots, \mathcal{F}_2(x), \sqrt{1-x^2}\mathcal{F}_{1}(x), \mathcal{F}_{0}(x)\}
\end{equation}
is a basis for $\Omega_{2m}$, we show that
\[
   \{\mathcal{F}_{2m+2}(x), \sqrt{1-x^2}\mathcal{F}_{2m+1}(x), \mathcal{F}_{2m}(x), \ldots, \mathcal{F}_2(x), \sqrt{1-x^2}\mathcal{F}_{1}(x), \mathcal{F}_{0}(x)\},
\]
is a basis for $\Omega_{2m+2}$.  Let $ c_0, c_1,c_2, \ldots, c_{2m+1}, c_{2m+2}$ be such that
\[
   c_{0}\,\mathcal{F}_{2m+2}(x) + c_{1}\,\sqrt{1-x^2}\mathcal{F}_{2m+1}(x) +  \ldots + c_{2m+1}\sqrt{1-x^2}\mathcal{F}_{1}(x) + c_{2m+2}\,\mathcal{F}_{0}(x) = 0.
\]
By considering the coefficients of $x^{2m+2}$ and $x^{2m+1} \sqrt{1-x^2}$, we have
\[
 \begin{array}{rl}
   c_0\,b_{2m+2,0}^{(0)} - c_1\,b_{2m+1,0}^{(1)} &= 0, \\[1ex]
   c_0\,b_{2m+2,0}^{(1)} + c_1\,b_{2m+1,0}^{(0)} &= 0.
 \end{array}
\]
The determinant of this system is $b_{2m+2,0}^{(0)}b_{2m+1,0}^{(0)}+b_{2m+2,0}^{(1)}b_{2m+1,0}^{(1)}$.  Since this determinant is different from zero,  we must have $c_0 = c_1 = 0$. This reduces our verification to finding $c_2, c_3, \ldots, c_{2m+2}$, such that
\[
   c_{2}\,\mathcal{W}_{2m}(x) + c_{3}\,\sqrt{1-x^2}\mathcal{W}_{2m-1}(x) +  \ldots + c_{2m+1}\sqrt{1-x^2}\mathcal{W}_{1}(x) + c_{2m+2}\,\mathcal{W}_{0}(x) = 0.
\]
Clearly $c_2=c_3 = \ldots = c_{2m+1}=c_{2m+2} = 0$, which follow from the assumption in (\ref{Assumption-1}). Thus, the results of the theorem for even indices follow by induction.  \hfill \qed

Now with the assumption $\gamma_0^2 > 0$, $\gamma_m^2 + \beta_m^2 > 0$, $m \geq 1$, we consider the functions $\mathcal{W}_m$ given by the recurrence formula (\ref{TTRR-Wm}).
It is easily seen that $\mathcal{W}_m$  takes the form
\begin{equation} \label{Decomposition-of-Wm}
   \mathcal{W}_m(x) = A_{m}^{(0)}(x) + \sqrt{1-x^2}  A_{m}^{(1)}(x),
\end{equation}
where $A_{m}^{(0)}$ and $A_{m}^{(1)}$ are, respectively, polynomials of degree at most $m$ and $m-1$, such that
\[
     A_{m}^{(0)}(-x) = (-1)^{m} A_{m}^{(0)}(x)\quad \mbox{and} \quad A_{m}^{(1)}(-x) = (-1)^{m-1} A_{m}^{(1)}(x).
\]
Setting
\begin{equation} \label{Coeffs-Decomposition-of-Wm}
  A_{m}^{(0)}(x) = \sum_{j=0}^{\lfloor m/2 \rfloor}  a_{m,\,2j}^{(0)}\,x^{m-2j} \ \ \mbox{and} \ \ A_{m}^{(1)}(x) = \sum_{j=0}^{\lfloor (m-1)/2 \rfloor}  a_{m,\,2j}^{(1)}\,x^{m-1-2j},
\end{equation}
we have the following.

\begin{theo} \label{Thm-Leading-Coeffs-Wm} Let $\gamma_0^2 > 0$, $\gamma_m^2 + \beta_m^2 > 0$, $m \geq 1$ in $(\ref{TTRR-Wm})$.  Then for the leading coefficients $a_{m,0}^{(0)}$ and $a_{m,0}^{(1)}$ and lead factor $\lambda_m =  (a_{m,0}^{(0)})^2 + (a_{m,0}^{(1)})^2$ of $\mathcal{W}_m$ the following hold.
 \[
   \left[
   \begin{array}{l}
     a_{m,0}^{(0)} \\[1ex]
     a_{m,0}^{(1)}
   \end{array}
   \right]
   =
   \left[
   \begin{array}{cc}
     \gamma_{m} & \beta_{m} \\[1ex]
     -\beta_{m} & \gamma_{m}
   \end{array}
   \right] \,
   \left[
   \begin{array}{l}
     a_{m-1,0}^{(0)} \\[1ex]
     a_{m-1,0}^{(1)}
   \end{array}
   \right],
   \quad m \geq 1,
 \]
with $a_{0,0}^{(0)} = 1$ and $a_{0,0}^{(1)} = 0$. Consequently, there hold
\[
 \begin{array}{ll}
     \lambda_0 = \gamma_0^2, \qquad \lambda_m = (\gamma_{m}^2+\beta_m^2)\lambda_{m-1}, \quad m \geq 1,
 \end{array}
\]
\[
 \begin{array}{l}
     \lambda_{m,1} = a_{m+1,0}^{(0)}\,a_{m,0}^{(0)} + a_{m+1,0}^{(1)}\,a_{m,0}^{(1)}  = \gamma_{m+1}\lambda_m,  \quad m \geq 1,
 \end{array}
\]
\[
 \begin{array}{l}
     a_{m+1,0}^{(0)}\,a_{m,0}^{(1)} - a_{m+1,0}^{(1)}\,a_{m,0}^{(0)}  = \beta_{m+1}\lambda_m, \quad m \geq 1,
 \end{array}
\]
\[
 \begin{array} l
   a_{m+1,0}^{(0)}\,a_{m-1,0}^{(0)} + a_{m+1,0}^{(1)}\,a_{m-1,0}^{(1)}
               = (\gamma_m\gamma_{m+1}-\beta_m\beta_{m+1}) \lambda_{m-1}, \quad m \geq 1
 \end{array}
\]
and
\[
 \begin{array} l
    a_{m+1,0}^{(0)}\,a_{m-1,0}^{(1)} - a_{m+1,0}^{(1)}\,a_{m-1,0}^{(0)}
    = (\gamma_{m+1}\beta_{m}+\gamma_m\beta_{m+1}) \lambda_{m-1}, \quad m \geq 1.
 \end{array}
\]
\end{theo}

\noindent{\bf Proof}.  From (\ref{TTRR-Wm}), (\ref{Decomposition-of-Wm}) and (\ref{Coeffs-Decomposition-of-Wm}), by equating the coefficients of $x^{m+1}$ and $x^{m}\sqrt{1-x^2}$,
\begin{equation} \label{Eq-SystemLeadingCoeffs}
 \begin{array}l
    a_{m+1,0}^{(0)} = \ \ \, \gamma_{m+1} a_{m,0}^{(0)} + \beta_{m+1}\,a_{m,0}^{(1)},  \\[2ex]
    a_{m+1,0}^{(1)} =  -\beta_{m+1}\,a_{m,0}^{(0)} + \gamma_{m+1}a_{m,0}^{(1)}.
 \end{array}
\end{equation}
From these the matrix formula in the theorem follows. The equalities and recurrence can be obtained as follows.

From the matrix formula,
\[
 \begin{array}{ll}
   (a_{m,0}^{(0)})^2 + (a_{m,0}^{(1)})^2  & = %
   \begin{array}{c}
     \left[ a_{m,0}^{(0)} \ \  a_{m,0}^{(1)}\right]\\[1ex]
        \quad
   \end{array} \!\!\!
   \left[
   \begin{array}{l}
     a_{m,0}^{(0)} \\[1ex]
     a_{m,0}^{(1)}
   \end{array}
   \right] \\[5ex]
   & = \begin{array}{c}
     \left[ a_{m-1,0}^{(0)} \ \  a_{m-1,0}^{(1)}\right]\\[1ex]
        \quad
   \end{array} \!\!\!
   \left[
   \begin{array}{cc}
     \gamma_{m} & -\beta_{m} \\[1ex]
     \beta_{m} & \gamma_{m}
   \end{array}
   \right]
   \left[
   \begin{array}{cc}
     \gamma_{m} & \beta_{m} \\[1ex]
     -\beta_{m} & \gamma_{m}
   \end{array}
   \right]
   \left[
   \begin{array}{l}
     a_{m-1,0}^{(0)} \\[1ex]
     a_{m-1,0}^{(1)}
   \end{array}
   \right].
 \end{array}
\]
Thus,
\[
 \begin{array}{ll}
   (a_{m,0}^{(0)})^2 + (a_{m,0}^{(1)})^2  & =  \begin{array}{c}
     \left[ a_{m-1,0}^{(0)} \ \  a_{m-1,0}^{(1)}\right]\\[1ex]
        \quad
   \end{array} \!\!\!
   \left[
   \begin{array}{cc}
     \gamma_{m}^2+\beta_m^2 & 0 \\[1ex]
     0 & \gamma_{m}^2+\beta_m^2
   \end{array}
   \right] \,
   \left[
   \begin{array}{l}
     a_{m-1,0}^{(0)} \\[1ex]
     a_{m-1,0}^{(1)}
   \end{array}
   \right] \\[5ex]
   & = (\gamma_{m}^2+\beta_m^2)\,\big[(a_{m-1,0}^{(0)})^2 + (a_{m-1,0}^{(1)})^2\big] = (\gamma_{m}^2+\beta_m^2)\,\lambda_{m-1},
 \end{array}
\]
which gives the recurrence formula for $\lambda_m$. Similarly,

\[
 \begin{array}{ll}
   a_{m+1,0}^{(0)}\,a_{m,0}^{(0)} + a_{m+1,0}^{(1)}\,a_{m,0}^{(1)}  & = %
   \begin{array}{c}
     \left[ a_{m,0}^{(0)} \ \  a_{m,0}^{(1)}\right]\\[1ex]
        \quad
   \end{array} \!\!\!
   \left[
   \begin{array}{l}
     a_{m+1,0}^{(0)} \\[1ex]
     a_{m+1,0}^{(1)}
   \end{array}
   \right] \\[5ex]
   & = \begin{array}{c}
     \left[ a_{m,0}^{(0)} \ \  a_{m,0}^{(1)}\right]\\[1ex]
        \quad
   \end{array} \!\!\!
   \left[
   \begin{array}{cc}
     \gamma_{m+1} & \beta_{m+1} \\[1ex]
     -\beta_{m+1} & \gamma_{m+1}
   \end{array}
   \right]
   \left[
   \begin{array}{l}
     a_{m,0}^{(0)} \\[1ex]
     a_{m,0}^{(1)}
   \end{array}
   \right] \\[5ex]
   & = \gamma_{m+1}\big[(a_{m,0}^{(0)})^2 + (a_{m,0}^{(1)})^2\big] = \gamma_{m+1}\lambda_{m}.
 \end{array}
\]
From (\ref{Eq-SystemLeadingCoeffs}),
\[
  \begin{array}l
   a_{m+1,0}^{(0)}\,a_{m,0}^{(1)} - a_{m+1,0}^{(1)}\,a_{m,0}^{(0)} \\[2ex]
   \qquad \qquad = [\gamma_{m+1}a_{m,0}^{(0)}+\beta_{m+1}a_{m,0}^{(1)}]a_{m,0}^{(1)} - [-\beta_{m+1}a_{m,0}^{(0)}+\gamma_{m+1}a_{m,0}^{(1)}]a_{m,0}^{(0)} \\[2ex]
   \qquad \qquad =  \beta_{m+1}\big[(a_{m,0}^{(0)})^2 +  (a_{m,0}^{(1)})^2\big] = \beta_{m+1} \lambda_m .
  \end{array}
\]
Again from (\ref{Eq-SystemLeadingCoeffs}),
\[
  \begin{array}l
   a_{m+1,0}^{(0)}\,a_{m-1,0}^{(0)} + a_{m+1,0}^{(1)}\,a_{m-1,0}^{(1)} \\[2ex]
   \qquad \qquad = [\gamma_{m+1} a_{m,0}^{(0)} +\beta_{m+1}a_{m,0}^{(1)}]a_{m-1,0}^{(0)} + [-\beta_{m+1}a_{m,0}^{(0)} + \gamma_{m+1} a_{m,0}^{(1)}]a_{m-1,0}^{(1)} \\[2ex]
   \qquad \qquad = \gamma_{m+1}\big[a_{m,0}^{(0)}\,a_{m-1,0}^{(0)} + a_{m,0}^{(1)}\,a_{m-1,0}^{(1)}\big] - \beta_{m+1}\big[a_{m,0}^{(0)}\,a_{m-1,0}^{(1)} - a_{m,0}^{(1)}\,a_{m-1,0}^{(0)}\big] \\[2ex]
   \qquad \qquad = (\gamma_m\gamma_{m+1}-\beta_m\beta_{m+1}) \lambda_{m-1}.
  \end{array}
\]
Similarly, the value associated with $a_{m+1,0}^{(0)}\,a_{m-1,0}^{(1)} - a_{m+1,0}^{(1)}\,a_{m-1,0}^{(0)}$ is also obtained.
\hfill \qed

Clearly, with the assumptions $\gamma_0^2 > 0$ and $\beta_m^2 + \gamma_m^2 > 0$, $m \geq 1$, there holds
\[
  \lambda_{m} = (a_{m,0}^{(0)})^2 + (a_{m,0}^{(1)})^2 > 0, \quad m \geq 0,
\]
which means the leading coefficients of $A_{m}^{(0)}$ and $A_{m}^{(1)}$ can not be zero simultaneously and $\mathcal{W}_m \in \Omega_m$ is of exact degree $m$.

With a more restrictive condition than $\gamma_0^2 > 0$ and $\beta_m^2 + \gamma_m^2 > 0$, $m \geq 1$, we  can state the following theorem.

\begin{theo} \label{Thm-2-Basis}
 Let $\gamma_m \neq 0$, $m \geq 0$ and let  $\{\mathcal{W}_m\}$ be the sequence of functions given by the recurrence formula $(\ref{TTRR-Wm})$. Then for any  $m \geq 0$, \\[-1ex]

\noindent 1.\ \ A basis for $\Omega_{2m}$ is
\[
   \{\mathcal{W}_{2m}(x), \sqrt{1-x^2}\mathcal{W}_{2m-1}(x), \mathcal{W}_{2m-2}(x), \ldots, \mathcal{W}_2(x), \sqrt{1-x^2}\mathcal{W}_{1}(x), \mathcal{W}_{0}(x)\};
\]

\noindent 2.\ \ A basis for $\Omega_{2m+1}$ is
\[
   \{\mathcal{W}_{2m+1}(x), \sqrt{1-x^2}\mathcal{W}_{2m}(x), \mathcal{W}_{2m-1}(x), \ldots, \mathcal{W}_3(x), \sqrt{1-x^2}\mathcal{W}_{2}(x), \mathcal{W}_{1}(x), \sqrt{1-x^2}\mathcal{W}_{0}(x)\}.
\]
\end{theo}

\noindent{\bf Proof}. From Theorem \ref{Thm-Leading-Coeffs-Wm} we observe that, with $\gamma_m \neq 0$, $m \geq 0$,  the leading coefficients of $\mathcal{W}_{m}$, $m \geq 0$ are such that $\lambda_{m,1} = a_{m+1,0}^{(0)}\,a_{m,0}^{(0)} + a_{m+1,0}^{(1)}\,a_{m,0}^{(1)} \neq 0$, $m \geq 0$. Hence, the  present theorem follows from Theorem \ref{Thm-Basis}.  \hfill \qed

Finally, by denoting the self-inversive polynomial associated with $\mathcal{W}_{m}(x)$ by  $K_m(z)$,  we have the following.
\begin{theo} \label{Thm-Leading-Coeff-SIP}
Let $ K_m(z) = \sum_{j=0}^{m} k_j^{(m)} z^j$ be such that
\[
   e^{-im\theta/2} K_m(e^{i\theta}) = \mathcal{W}_{m}(x),
\]
where $x = \cos(\theta/2)$. Then
\[
     k_0^{(m)} = \overline{k_m^{(m)}} = 2^{-m}[a_{m,0}^{(0)} + i\, a_{m,0}^{(1)}] \quad \mbox{and} \quad |k_0^{(m)}|^2 = |k_m^{(m)}|^2 = 2^{-2m}\lambda_m.
\]
\end{theo}

\setcounter{equation}{0}
\section{Orthogonal properties associated with $\mathcal{W}_m$  }

Let  $\psi$ be a positive measure on $[-1,1]$. We consider the sequence of functions $\{\mathcal{W}_m\}$, where $\mathcal{W}_m \in \Omega_m$ is of exact degree $m$, is such that
\begin{equation} \label{Orthogonality-for-Wm}
  \begin{array}{l}
    \dsp \int_{-1}^{1} \mathcal{W}_{2n}(x)\,\mathcal{W}_{2m}(x)\, \sqrt{1-x^2}\,d\psi(x) = \rho_{2m}\,\delta_{n,m} ,  \\[2ex]
    \dsp \int_{-1}^{1} \mathcal{W}_{2n+1}(x)\,\mathcal{W}_{2m+1}(x)\, \sqrt{1-x^2}\,d\psi(x) =  \rho_{2m+1} \,\delta_{n,m},\\[2ex]
    \dsp \int_{-1}^{1} \mathcal{W}_{2n+1}(x)\,\mathcal{W}_{2m}(x)\, d\psi(x) = 0,
  \end{array}
\end{equation}
for $n,m = 0, 1, 2, \ldots \ $ .

We use the notation  $\mathcal{W}_0(x) = a_{0,0}^{(0)}$,
\[
  \mathcal{W}_m(x)  = \sum_{j=0}^{\lfloor m/2 \rfloor}  a_{m,\,2j}^{(0)}\,x^{m-2j} \ \ + \sqrt{1-x^2} \sum_{j=0}^{\lfloor (m-1)/2 \rfloor}  a_{m,\,2j}^{(1)}\,x^{m-1-2j}, \quad m \geq 1,
\]
and
\[ \lambda_m = (a_{m,0}^{(0)})^2 + (a_{m,0}^{(1)})^2 \quad \mbox{and} \quad \lambda_{m,1} =
 a_{m+1,0}^{(0)}a_{m,0}^{(0)}+a_{m+1,0}^{(1)}a_{m,0}^{(1)}, \quad m \geq 0,
\]
with $a_{0,0}^{(1)} = 0$.

Observe that, in the case $a_{m, 2j}^{(1)} = 0$, $j=0, 1,  \ldots, \lfloor (m-1)/2\rfloor$, $m \geq 1$, then $\{\mathcal{W}_m\}$ are symmetric polynomials and (\ref{Orthogonality-for-Wm}) reduces to the orthogonality of symmetric polynomials with respect to $\sqrt{1-x^2}\,d\psi(x)$ in $[-1, 1]$. For references to some of the classical texts on orthogonal polynomials on the real line we cite \cite{Chihara-book, Ismail-book-2005, Szego-book}.

\begin{theo} \label{Thm-Orthogonality-Wm}
The sequence of functions $\{\mathcal{W}_m\}$ that satisfies the orthogonality property $(\ref{Orthogonality-for-Wm})$, where $\mathcal{W}_m \in \Omega_m$ and is of exact degree $m$, exists and satisfies the three term recurrence formula
\begin{equation} \label{TTRR-Wm-thm}
  \begin{array}{l}
  \mathcal{W}_0(x) = \gamma_0,\qquad \mathcal{W}_1(x) = \big[\gamma_1x - \beta_{1} \sqrt{1-x^2}\big]\gamma_0,  \\[2ex]
  \mathcal{W}_{m+1}(x) = \big[\gamma_{m+1}x - \beta_{m+1} \sqrt{1-x^2}\big] \mathcal{W}_{m}(x) - \alpha_{m+1} \mathcal{W}_{m-1}(x), \quad m \geq 1,
  \end{array}
\end{equation}
where \ \ $\gamma_m \neq 0$, $m \geq 0$, \ \ $\beta_1 = \gamma_1 \,\rho_{0}^{-1} \int_{-1}^{1} x\, \mathcal{W}_{0}^2(x)\, d\psi(x)$,
\begin{equation} \label{Coef-TTRR-Wm-thm}
  \begin{array}l
     \dsp \beta_{m+1} = \gamma_{m+1}\,
     \frac{1}{\rho_{m}} \int_{-1}^{1} x\, \mathcal{W}_{m}^2(x)\, d\psi(x), \quad m \geq 1  \qquad \mbox{and}
      \\[3ex]
     \dsp \alpha_{m+1} = \frac{1}{\rho_{m-1}}\int_{-1}^{1} \big[\gamma_{m+1}x - \beta_{m+1}\sqrt{1-x^2}\big]\, \mathcal{W}_{m-1}(x)\, \mathcal{W}_{m}(x)\sqrt{1-x^2}\, d\psi(x), \quad m \geq 1.
  \end{array}
\end{equation}
Here, $\rho_m = \int_{-1}^{1}  \mathcal{W}_{m}^2(x)\,\sqrt{1-x^2}\, d\psi(x)$, $m \geq 0$.
\end{theo}

\noindent{\bf Proof}.  It follows from Theorem  \ref{Thm-Leading-Coeffs-Wm} that, with $\gamma_m \neq 0$, $m \geq 0$, the function $\mathcal{W}_{m}$  obtained from the three term recurrence formula (\ref{TTRR-Wm-thm}) is of exact degree $m$. Hence, to prove the existence part of the above theorem, we show that the sequence of functions $\{\mathcal{W}_{m}\}$ generated by the above three term recurrence formula satisfies $(\ref{Orthogonality-for-Wm})$.

Since $\beta_1 = \gamma_1 \,\rho_{0}^{-1} \int_{-1}^{1} x\, \mathcal{W}_{0}^2(x)\, d\psi(x)$, with $\mathcal{W}_1(x) = \big[\gamma_1x - \beta_{1} \sqrt{1-x^2}\big]\gamma_0$ there follows $\int_{-1}^{1} \mathcal{W}_{0}(x)\,\mathcal{W}_{1}(x)\, d\psi(x) = 0$.

Since $\beta_{2}$ and $\alpha_{2}$ are as in (\ref{Coef-TTRR-Wm-thm}), with
\[
   \mathcal{W}_{2}(x) = \big[\gamma_2 x - \beta_{2} \sqrt{1-x^2}\big] \mathcal{W}_{1}(x) - \alpha_{2} \mathcal{W}_{0}(x),
\]
we have
\[
  \int_{-1}^{1} \mathcal{W}_{1}(x)\,\mathcal{W}_{2}(x)\, d\psi(x) = 0 \quad \mbox{and} \quad \int_{-1}^{1} \mathcal{W}_{0}(x)\,\mathcal{W}_{2}(x)\, \sqrt{1-x^2}\,d\psi(x) = 0.
\]

Now suppose that for $N \geq 2$ the sequence of functions  $\{\mathcal{W}_m\}_{m=0}^{N}$ obtained from the three  term recurrence formula (\ref{TTRR-Wm-thm}) satisfies  $(\ref{Orthogonality-for-Wm})$.

Since \ $x\mathcal{W}_{N-2k}(x) \in \Omega_{N+1-2k}$, by Theorem \ref{Thm-2-Basis} there exists $c_0, c_1, \ldots, c_{N-2k}$ such  that
\[
    x \mathcal{W}_{N-2k}(x) = c_0\mathcal{W}_{N+1-2k}(x)+ c_1\sqrt{1-x^2}\mathcal{W}_{N-2k}(x)+ c_2\mathcal{W}_{N-1-2k}(x) + \ldots \ .
\]
Hence, we have
\begin{equation*}  \label{Orthogonality-1-for-Wm}
     \int_{-1}^{1} x  \mathcal{W}_{N-2k}(x)\, \mathcal{W}_{N}(x)\, d\psi(x) = 0, \quad k=1, 2, \ldots \lfloor N/2 \rfloor.
\end{equation*}
Likewise, since \ $x \sqrt{1-x^2} \mathcal{W}_{N-1-2k}(x) \in \Omega_{N+1-2k}$, we have
\begin{equation}  \label{Orthogonality-2-for-Wm}
     \int_{-1}^{1} x \sqrt{1-x^2} \mathcal{W}_{N-1-2k}(x)\, \mathcal{W}_{N}(x)\, d\psi(x) = 0, \quad k=1, 2, \ldots \lfloor (N-1)/2 \rfloor \ \
\end{equation}
and, since $(1-x^2) \mathcal{W}_{N-1-2k}(x) \in \Omega_{N+1-2k}$, we also have
\begin{equation}  \label{Orthogonality-3-for-Wm}
     \int_{-1}^{1} (1-x^2) \mathcal{W}_{N-1-2k}(x)\, \mathcal{W}_{N}(x)\, d\psi(x) = 0, \quad k=1, 2, \ldots \lfloor (N-1)/2 \rfloor.
\end{equation}
Note that in (\ref{Orthogonality-2-for-Wm}) and (\ref{Orthogonality-3-for-Wm}) the value of $N$ is assumed to be $\geq 3$.

Hence, from
\[
   \mathcal{W}_{N+1}(x) = \big[\gamma_{N+1}x - \beta_{N+1} \sqrt{1-x^2}\big] \mathcal{W}_{N}(x) - \alpha_{N+1} \mathcal{W}_{N-1}(x),
\]
it follows that
\[
 \begin{array}{l}
  \int_{-1}^{1} \mathcal{W}_{N-2k}(x)\,\mathcal{W}_{N+1}(x)\, d\psi(x) = 0, \quad k=1,2, \ldots, \lfloor N/2 \rfloor, \\[3ex]
  \int_{-1}^{1} \mathcal{W}_{N-1-2k}(x)\,\mathcal{W}_{N+1}(x)\, \sqrt{1-x^2} d\psi(x) = 0, \quad k=1,2, \ldots, \lfloor (N-1)/2 \rfloor.
 \end{array}
\]
Moreover, since $\alpha_{N+1}$ and $\beta_{N+1}$ are as in (\ref{Coef-TTRR-Wm-thm}), we also have
\[
 \begin{array}{l}
  \int_{-1}^{1} \mathcal{W}_{N}(x)\,\mathcal{W}_{N+1}(x)\, d\psi(x) = 0  \quad \mbox{and} \quad  \int_{-1}^{1} \mathcal{W}_{N-1}(x)\,\mathcal{W}_{N+1}(x)\, \sqrt{1-x^2} d\psi(x) = 0.
 \end{array}
\]
Thus, by induction we  conclude that the sequence of functions given by the three term recurrence formula satisfies $(\ref{Orthogonality-for-Wm})$.

On the other hand, to show that any  sequence of functions $\{\mathcal{W}_m\}$, where $\mathcal{W}_m \in \Omega_m$ is of exact degree $m$, for which  (\ref{Orthogonality-for-Wm}) holds, must also satisfy the three term recurrence formula (\ref{TTRR-Wm-thm}), we proceed as follows.

Clearly we can write $\mathcal{W}_0(x) = \gamma_0 \neq 0$ and $\mathcal{W}_1(x) = a_{1,0}^{(0)}x + a_{1,0}^{(1)}\sqrt{1-x^2} =  [\gamma_1x - \beta_1 \sqrt{1-x^2}]\gamma_0$. Then for $\int_{-1}^{1} \mathcal{W}_{0}(x)\,\mathcal{W}_{1}(x)\, d\psi(x) = 0$ to hold such that $\mathcal{W}_1(x)$ is of exact degree 1, one must have $\gamma_1 \neq 0$ and $\beta_1 = \gamma_1 \,\rho_{0}^{-1} \int_{-1}^{1} x\, \mathcal{W}_{0}^2(x)\, d\psi(x)$.

With the next element $\mathcal{W}_2$ of the given orthogonality sequence, let $\gamma_{2}$ and $\beta_{2}$ be such that \linebreak $\mathcal{W}_{2}(x) - \big[\gamma_2x - \beta_{2}\sqrt{1-x^2}\big]\mathcal{W}_{1}(x) \in \Omega_{0}$. With respect to the  leading coefficients of $\mathcal{W}_{2}$ and $\mathcal{W}_{1}$, the elements  $\gamma_{2}$ and $\beta_{2}$  must satisfy
\[
   \left[
   \begin{array}{l}
     a_{2,0}^{(0)} \\[1ex]
     a_{2,0}^{(1)}
   \end{array}
   \right]
   =
   \left[
   \begin{array}{cc}
     a_{1,0}^{(0)} & a_{1,0}^{(1)} \\[1ex]
     a_{1,0}^{(1)} & -a_{1,0}^{(0)}
   \end{array}
   \right] \,
   \left[
   \begin{array}{l}
     \gamma_{2} \\[2ex]
     \beta_{2}
   \end{array}
   \right].
\]
The determinant of this system is $-[(a_{1,0}^{(0)})^2+ (a_{1,0}^{(1)})^2]$. Since $\mathcal{W}_{1}$ is of exact degree $1$ this determinant is different from zero. Hence, the values of  $\gamma_{2}$ and $\beta_{2}$ are uniquely found.  Writing
\[
   \mathcal{W}_{2}(x) = \big[\gamma_2x - \beta_{2}\sqrt{1-x^2}\big]\mathcal{W}_{1}(x) - \alpha_2 \mathcal{W}_{0}(x),
\]
we find, with the orthogonality and with the additional observation that $\mathcal{W}_{2}$ is of exact degree 2, that $\gamma_2 \neq 0$,
\begin{equation*}
  \begin{array}l
     \dsp \beta_{2} = \gamma_{2}\,
     \frac{1}{\rho_{1}} \int_{-1}^{1} x\, \mathcal{W}_{1}^2(x)\, d\psi(x)  \qquad \mbox{and}
      \\[3ex]
     \dsp \alpha_{2} = \frac{1}{\rho_{0}}\int_{-1}^{1} \big[\gamma_{2}x - \beta_{2}\sqrt{1-x^2}\big]\, \mathcal{W}_{0}(x)\, \mathcal{W}_{1}(x)\sqrt{1-x^2}\, d\psi(x).
  \end{array}
\end{equation*}

Now for  $N \geq 2$ assume that the orthogonal functions $\mathcal{W}_m$, $m=0,1, \ldots, N$ satisfy  the three term recurrence formula (\ref{TTRR-Wm-thm}), with $\gamma_m \neq 0$, $m=0,1 , \ldots, N$.  Let $\gamma_{N+1}$ and $\beta_{N+1}$ be such that $\mathcal{W}_{N+1}(x) - \big[\gamma_{N+1} x - \beta_{N+1}\sqrt{1-x^2}\big]\mathcal{W}_{N}(x) \in \Omega_{N-1}$. With respect to the  leading coefficients of $\mathcal{W}_{N+1}$ and $\mathcal{W}_{N}$, the elements  $\gamma_{N+1}$ and $\beta_{N+1}$  must satisfy
\[
   \left[
   \begin{array}{l}
     a_{N+1,0}^{(0)} \\[1ex]
     a_{N+1,0}^{(1)}
   \end{array}
   \right]
   =
   \left[
   \begin{array}{cc}
     a_{N,0}^{(0)} & a_{N,0}^{(1)} \\[1ex]
     a_{N,0}^{(1)} & -a_{N,0}^{(0)}
   \end{array}
   \right] \,
   \left[
   \begin{array}{l}
     \gamma_{N+1} \\[1ex]
     \beta_{N+1}
   \end{array}
   \right].
\]
The determinant of this system is $-[(a_{N,0}^{(0)})^2+ (a_{N,0}^{(1)})^2]$, which is different form zero because $\mathcal{W}_{N}$ is of exact degree $N$. Hence, the values for $\gamma_{N+1}$ and $\beta_{N+1}$ are uniquely found.

Now using Theorem \ref{Thm-2-Basis}, there exist $c_0, c_1, \ldots, c_{N-1}$ such that
\[
 \begin{array} {ll}
    \mathcal{W}_{N+1}(x) = & \big[\gamma_{N+1} x - \beta_{N+1}\sqrt{1-x^2}\big]\mathcal{W}_{N}(x) \\[2ex]
    & \qquad \quad + \ c_0 \mathcal{W}_{N-1}(x) + c_1 \sqrt{1-x^2} \mathcal{W}_{N-2}(x) + c_2 \mathcal{W}_{N-3}(x) + \ldots \ .
 \end{array}
\]
Applications of the orthogonality properties (\ref{Orthogonality-for-Wm}) of the sequence $\{\mathcal{W}_{m}\}_{m=0}^{N+1}$, together with  the observation that $\mathcal{W}_{N+1}$ is of exact degree $N+1$, lead to  $c_1 = c_2 = \ldots = c_{N-1} = 0$,  $\gamma_{N+1} \neq 0$,
\begin{equation*}
  \begin{array}l
     \dsp \beta_{N+1} = \gamma_{N+1}\,
     \frac{1}{\rho_{N}} \int_{-1}^{1} x\, \mathcal{W}_{N}^2(x)\, d\psi(x)  \qquad \mbox{and}
      \\[3ex]
     \dsp \alpha_{N+1} = -c_0 = \frac{1}{\rho_{N-1}}\int_{-1}^{1} \big[\gamma_{N+1}x - \beta_{N+1}\sqrt{1-x^2}\big]\, \mathcal{W}_{N-1}(x)\, \mathcal{W}_{N}(x)\sqrt{1-x^2}\, d\psi(x).
  \end{array}
\end{equation*}
This concludes the proof of the Theorem.  \hfill \qed

Following as in Theorem \ref{Thm-Leading-Coeffs-Wm} we have $\gamma_{m+1} = \lambda_{m,1}/\lambda_m$. Thus, the orthogonal functions $\mathcal{W}_{m}$, $m \geq 1$, are such that
\[
    \lambda_{m,1} = a_{m+1,0}^{(0)}\,a_{m,0}^{(0)} + a_{m+1,0}^{(1)}\,a_{m,0}^{(1)} \neq 0, \quad m \geq 0.
\]

\begin{coro} \label{Coro-Orthogonality-Wm}

The sequence of functions $\{\mathcal{W}_m\}$, where $\mathcal{W}_m \in \Omega_m$ and is of exact degree $m$, satisfies the orthogonality property $(\ref{Orthogonality-for-Wm})$ if and only if,\ for $m\geq 1$,
\begin{equation} \label{Eq-Coro-Orthogonality-for-Wm}
 \begin{array}l
  \dsp \int_{-1}^{1} \mathcal{F}(x)\,\mathcal{W}_{m}(x)\, d\psi(x) = 0 \quad \mbox{whenever} \quad \mathcal{F} \in \Omega_{m-1}.
 \end{array}
\end{equation}

\end{coro}

\noindent {\bf Proof}. First we assume that (\ref{Eq-Coro-Orthogonality-for-Wm}) holds. Observe that if $\mathcal{F} \in \Omega_{m+1-2k}$ for $k=1,2, \ldots, \lfloor (m+1)/2 \rfloor$ then $\mathcal{F} \in \Omega_{m-1}$.  Hence, $\mathcal{W}_{m+1-2k}(x) \in \Omega_{m+1-2k}$ and $\mathcal{W}_{m-2k}(x) \sqrt{1-x^2} \in \Omega_{m+1-2k}$ leads immediately to $(\ref{Orthogonality-for-Wm})$.

On the other hand, if $(\ref{Orthogonality-for-Wm})$ holds then from Theorem \ref{Thm-Orthogonality-Wm} and from  Theorem \ref{Thm-2-Basis} we can write
\[
   \mathcal{F}_{m+1-2k}(x) = c_0 \mathcal{W}_{m+1-2k}(x) + c_1 \sqrt{1-x^2} \mathcal{W}_{m-2k}(x) + c_2 \mathcal{W}_{m-1-2k}(x) + \ldots \ .
\]
for $k=1,2, \ldots, \lfloor (m+1)/2 \rfloor$.  Hence, (\ref{Eq-Coro-Orthogonality-for-Wm}) is immediate.               \hfill \qed

It is quite straight forward that another way to present the above corollary is the following.

\begin{coro} \label{Coro2-Orthogonality-Wm}

The sequence of functions $\{\mathcal{W}_m\}$, where $\mathcal{W}_m \in \Omega_m$ and is of exact degree $m$, satisfies the orthogonality property $(\ref{Orthogonality-for-Wm})$ if and only if, \ for $m\geq 1$,
\begin{equation*} \label{Eq-Coro2-Orthogonality-for-Wm}
 \begin{array}l
  \dsp \int_{-1}^{1} B^{(0)}(x)\,\mathcal{W}_{m}(x)\, d\psi(x) = 0 \quad \mbox{and} \quad \int_{-1}^{1} B^{(1)}(x)\,\mathcal{W}_{m}(x)\, \sqrt{1-x^2}\,d\psi(x) = 0,
 \end{array} \quad
\end{equation*}
where $B^{(0)} \in \mathbb{P}_{m-1}$ and $B^{(1)} \in \mathbb{P}_{m-2}$ satisfy
\[
  \begin{array}l
    B^{(0)}(-x) = (-1)^{m-1}B^{(0)}(x) \ \ \ \mbox{and}  \ \ \ B^{(1)}(-x) = (-1)^{m-2}B^{(1)}(x).
  \end{array}
\]
\end{coro}

The following corollary provides one other way to express the orthogonality (\ref{Orthogonality-for-Wm}) of the sequence  $\{\mathcal{W}_n\}$.

\begin{coro} \label{Coro3-Orthogonality-Wm}

The sequence of functions $\{\mathcal{W}_m\}$, where $\mathcal{W}_m \in \Omega_m$ and is of exact degree $m$, satisfies the orthogonality property $(\ref{Orthogonality-for-Wm})$ if and only if,  for \ $m\geq 1$,
\begin{equation} \label{Eq-Coro3-Orthogonality-for-Wm}
 \begin{array}l
  \dsp \int_{-1}^{1} (x + i \sqrt{1-x^2}\,)^{-m+1+2s}\,\mathcal{W}_{m}(x)\, d\psi(x) = 0, \quad s=0,1,\ldots, m-1.
 \end{array} \quad
\end{equation}
\end{coro}

\noindent {\bf Proof}. We just give the proof for $m = 2n$ and the proof for $m = 2n+1$ is similar.

Since $(x + i \sqrt{1-x^2}\,)(x - i \sqrt{1-x^2}\,) = 1$ the orthogonality (\ref{Eq-Coro3-Orthogonality-for-Wm}) for $m =2n$ can be written as

\[
 \begin{array}l
  \dsp \int_{-1}^{1} (x \pm i \sqrt{1-x^2}\,)^{2l+1}\,\mathcal{W}_{2n}(x)\, d\psi(x) = 0, \quad l=0,1,\ldots, n-1.
 \end{array} \quad
\]
Observe that
\[
  \begin{array}l
    (x \pm i \sqrt{1-x^2}\,)^{2l+1} \\[1ex]
   \qquad \qquad \quad \dsp = \sum_{k=0}^{l} \left(\!\!\!\begin{array}l 2l+1 \\ 2k+1 \end{array}\!\!\!\right) x^{2k+1} (x^2-1)^{l-k} \pm i \sqrt{1-x^2}\, \sum_{k=0}^{l} \left(\!\!\!\begin{array}l  2l+1 \\ \ \ 2k \end{array}\!\!\!\right) x^{2k} (x^2-1)^{l-k}.
  \end{array}
\]
Since the polynomials represented by the above sums are, respectively, odd and even polynomials of exact degrees $2l+1$ and $2l$, the required result follows from Corollary \ref{Coro2-Orthogonality-Wm}.  \hfill \qed

The following theorem, in addition to showing some further orthogonality properties of the functions $\{\mathcal{W}_m\}$  given by Theorem \ref{Thm-Orthogonality-Wm}, also gives another expression for  the coefficients $\alpha_n$ given in Theorem \ref{Thm-Orthogonality-Wm}.

\begin{theo} \label{Thm-FurtherOrthogonality-Wm}
\[
  \int_{-1}^{1} x\sqrt{1-x^2}\, \mathcal{W}_{1}(x)\,\mathcal{W}_{2}(x)\,d\psi(x) = \frac{\gamma_{2}\gamma_{3} - \beta_{2}\beta_{3}} {(\gamma_2^2+\beta_{2}^2)\gamma_3}\,\rho_{2},
\]
\[
  \int_{-1}^{1} (1-x^2)\, \mathcal{W}_{1}(x)\,\mathcal{W}_{2}(x)\, d\psi(x) = -\frac{\gamma_{3}\beta_{2} + \gamma_{2}\beta_{3}} {(\gamma_{2}^2+\beta_{2}^2)\gamma_{3}}\,\rho_{2}
\]
and for $m \geq 2$,
\[
  \begin{array}l
    \dsp \int_{-1}^{1} x\sqrt{1-x^2}\, \mathcal{W}_{m-1-2k}(x)\,\mathcal{W}_{m}(x)\,d\psi(x) = \left\{ \begin{array}{cl}
       \dsp \frac{\gamma_{m}\gamma_{m+1}-\beta_{m}\beta_{m+1}}{(\gamma_{m}^2+\beta_{m}^2)\gamma_{m+1}}\,\rho_{m}, & k=0, \\[3ex]
       \dsp 0, & 1 \leq k \leq \lfloor(m-1)/2\rfloor,
    \end{array}
    \right. \\[5ex]
  \end{array}
\]
\[
  \begin{array}l
    \dsp \int_{-1}^{1} (1-x^2)\, \mathcal{W}_{m-1-2k}(x)\,\mathcal{W}_{m}(x)\, d\psi(x) = \left\{ \begin{array}{cl}
       \dsp -\frac{\gamma_{m+1}\beta_{m}+\gamma_{m}\beta_{m+1}}{(\gamma_{m}^2+\beta_{m}^2)\gamma_{m+1}}\,\rho_{m}, & k=0, \\[3ex]
       \dsp 0, & 1 \leq k \leq \lfloor(m-1)/2\rfloor,
    \end{array}
    \right. \\[5ex]
  \end{array}
\]
Consequently,
\[
    \alpha_{m+1} = \frac{\gamma_{m+1}^2+\beta_{m+1}^2}{\gamma_{m}^2+\beta_{m}^2}\, \frac{\gamma_{m}}{\gamma_{m+1}}\, \frac{\rho_{m}}{\rho_{m-1}}, \quad m \geq 1.
\]
\end{theo}

\noindent{\bf Proof}. Since $x\sqrt{1-x^2}\, \mathcal{W}_{m-1-2k}(x)$ and $(1-x^2)\, \mathcal{W}_{m-1-2k}(x)$ are in $\Omega_{m+1-2k}$,  the results corresponding to $1 \leq k \leq \lfloor(m-1)/2\rfloor$ follows from  Corollary \ref{Coro-Orthogonality-Wm}.

We now prove the result associated with  $\int_{-1}^{1} x\sqrt{1-x^2}\, \mathcal{W}_{m-1}(x)\,\mathcal{W}_{m}(x)\, d\psi(x)$.

Since $x\sqrt{1-x^2}\, \mathcal{W}_{m-1}(x) \in \Omega_{m+1}$, there exist $c_0, c_1, \ldots, c_{m+1}$ such that
\[
    x\sqrt{1-x^2}\, \mathcal{W}_{m-1}(x) = c_0 \mathcal{W}_{m+1}(x) + c_1 \sqrt{1-x^2} \mathcal{W}_{m}(x) + c_2 \mathcal{W}_{m-1}(x) + \ldots \
\]
and that $\int_{-1}^{1} x\sqrt{1-x^2}\, \mathcal{W}_{m-1}(x)\,\mathcal{W}_{m}(x)\, d\psi(x) = c_1 \rho_{m}$.  Moreover, comparing the leading coefficients on both sides
\[
 \begin{array}{rl}
  -a_{m-1,0}^{(1)} & = c_0\,a_{m+1,0}^{(0)} - c_1\,a_{m,0}^{(1)},  \\[1ex]
   a_{m-1,0}^{(0)} & = c_0\,a_{m+1,0}^{(1)} + c_1\,a_{m,0}^{(0)} .
 \end{array}
\]
 This gives $c_1 = \big[a_{m+1,0}^{(0)}a_{m-1,0}^{(0)}+a_{m+1,0}^{(1)}a_{m-1,0}^{(1)} \big]/\big[a_{m+1,0}^{(0)}a_{m,0}^{(0)}+a_{m+1,0}^{(1)}a_{m,0}^{(1)}\big]$. Thus, from Theorem \ref{Thm-Leading-Coeffs-Wm} we have
\[
    c_1 = \frac{\gamma_{m}\gamma_{m+1}-\beta_{m}\beta_{m+1}}{(\gamma_{m}^2+\beta_{m}^2)\gamma_{m+1}}
\]
and the required result for $\int_{-1}^{1} x\sqrt{1-x^2}\, \mathcal{W}_{m-1}(x)\,\mathcal{W}_{m}(x)\, d\psi(x)$. 

The result associated with $\int_{-1}^{1} (1-x^2)\, \mathcal{W}_{m-1}(x)\,\mathcal{W}_{m}(x)\, d\psi(x)$ is obtained similarly, and the latter result of the theorem  then follows from (\ref{Coef-TTRR-Wm-thm}). \hfill \qed

Now we can state the following theorem with respect to the zeros of $\mathcal{W}_m(x)$.

\begin{theo} \label{Thm-Zeros-Wm}
Let $\{\mathcal{W}_m(x)\}$ be the sequence of functions defined as in Theorem $\ref{Thm-Orthogonality-Wm}$.  Then for $m \geq 1$, the function  $\mathcal{W}_m(x)$ has exactly $m$ distinct zeros in $(-1,1)$.
\end{theo}

\noindent {\bf Proof}. From (\ref{Orthogonality-for-Wm}) since at least one of the integrals
\[
 \int_{-1}^{1} \mathcal{W}_0(x)\mathcal{W}_{m}(x)d\psi(x) \quad \mbox{and} \quad \int_{-1}^{1} \mathcal{W}_0(x)\mathcal{W}_{m}(x)\sqrt{1-x^2}d\psi(x)
\]
is zero, we can say that $\mathcal{W}_{m}(x)$ changes sign at least once in $(-1,1)$. According to Theorem \ref{Thm-Max-Zeros} the number of sign changes of $\mathcal{W}_m(x)$ also can not exceed $m$.

Suppose that $\mathcal{W}_{m}(x)$ changes sign $k$ ($1 \leq k \leq m$)  times in $(-1,1)$, namely at the points $y_1, y_2, \ldots, y_k$. Let $\theta_j = 2 \arccos(y_j)$, $j=1,2, \ldots, k$ and consider the self-inversive polynomial (of degree $k$) defined by
\[
     q_k(z) =  e^{-i(k\pi+\theta_1 + \theta_2 +\ldots+\theta_k)/2} (z - e^{i\theta_1})(z - e^{i\theta_2}) \cdots (z - e^{i\theta_k}).
\]
Now, with $x = \cos(\theta/2)= (z^{1/2}+z^{-1/2})/2$, if we consider the function $\mathcal{F}_k(x)$ given by
\[ \mathcal{F}_k(x) = e^{-ik\theta/2}q_k(e^{i\theta}), \]
then by Lemma \ref{Lemma1} $\mathcal{F}_k(x) \in \Omega_k$ and further $\mathcal{F}_k(x)$ has exactly the $k$ zeros $y_1, y_2, \ldots, y_k$, which are the points of sign changes in $\mathcal{W}_m(x)$. Hence, the function $\mathcal{F}_k(x)\mathcal{W}_m(x)$ does not change sign in $(-1,1)$ which then leads to the conclusion  %
\[
   \int_{-1}^{1} \mathcal{F}_k(x)\mathcal{W}_{m}(x)d\psi(x) \neq 0 \quad \mbox{and} \quad \int_{-1}^{1} \mathcal{F}_k(x)\mathcal{W}_{m}(x)\sqrt{1-x^2}d\psi(x) \neq 0.
\]
On the other hand, if $k < m$ then from Corollary \ref{Coro-Orthogonality-Wm} at least one of the integrals
\[
   \int_{-1}^{1} \mathcal{F}_k(x)\mathcal{W}_{m}(x)d\psi(x)  \quad \mbox{and} \quad \int_{-1}^{1} \mathcal{F}_k(x)\mathcal{W}_{m}(x)\sqrt{1-x^2}d\psi(x)
\]
must be equal to zero, contradicting the earlier conclusion. Hence, the only possibility is $k =m$ and that $\mathcal{W}_m$ has $m$ sign changes in $(-1,1)$. That is, $\mathcal{W}_m$ has exactly $m$ zeros in $(-1,1)$. \hfill \qed

Now we look at the sequence of functions $\{\widehat{\mathcal{W}}_{m}\}$, obtained from the sequence of orthogonal functions  $\{\mathcal{W}_m\}$ by the scaling
\[
   \widehat{\mathcal{W}}_{0}(x) = \frac{1}{\gamma_0}\,\mathcal{W}_0(x), \qquad \widehat{\mathcal{W}}_{m}(x) = \frac{1}{\gamma_0\cdots\gamma_m}\,\mathcal{W}_m(x), \quad m \geq 0.
\]
By considering the properties of $\{\widehat{\mathcal{W}}_{m}\}$ we can state the following.

\begin{theo} \label{Thm-Orthogonality-Normalized-Wm}
Let the sequence of  functions $\{\widehat{\mathcal{W}}_{m}\}$ be given by
\begin{equation*} \label{TTRR-MWm-thm}
  \begin{array}{l}
  \widehat{\mathcal{W}}_0(x) = 1,\qquad \widehat{\mathcal{W}}_1(x) = x - \widehat{\beta}_{1} \sqrt{1-x^2},  \\[2ex]
  \widehat{\mathcal{W}}_{m+1}(x) = \big[x - \widehat{\beta}_{m+1} \sqrt{1-x^2}\big] \widehat{\mathcal{W}}_{m}(x) - \widehat{\alpha}_{m+1} \widehat{\mathcal{W}}_{m-1}(x), \quad m \geq 1,
  \end{array}
\end{equation*}
where \ \  $\dsp \widehat{\beta}_1 = \, \frac{1}{\widehat{\rho}_{0}} \int_{-1}^{1} x\, \widehat{\mathcal{W}}_{0}^2(x)\, d\psi(x)$ \ and for \ $m \geq 1$,
\begin{equation*} \label{Coef-TTRR-MWm-thm}
  \begin{array}l
     \dsp \widehat{\beta}_{m+1} =
     \frac{1}{\widehat{\rho}_{m}} \int_{-1}^{1} x\, \widehat{\mathcal{W}}_{m}^2(x)\, d\psi(x) \qquad \mbox{and}
      \\[3ex]
     \dsp \widehat{\alpha}_{m+1} =  \frac{1}{\widehat{\rho}_{m-1}}\int_{-1}^{1} \big[x - \widehat{\beta}_{m+1}\sqrt{1-x^2}\big]\, \widehat{\mathcal{W}}_{m-1}(x)\, \widehat{\mathcal{W}}_{m}(x)\sqrt{1-x^2}\, d\psi(x) = \frac{1+\widehat{\beta}_{m+1}^2}{1+\widehat{\beta}_{m}^2}\, \frac{\widehat{\rho}_{m}}{\widehat{\rho}_{m-1}}.
  \end{array}
\end{equation*}
Here, $\widehat{\rho}_m = \int_{-1}^{1}  \widehat{\mathcal{W}}_{m}^2(x)\,\sqrt{1-x^2}\, d\psi(x)$, $m \geq 0$.
Then $\{\widehat{\mathcal{W}}_{m}\}$, where $\widehat{\mathcal{W}}_{m} \in \Omega_m$ and is of exact degree $m$,  satisfies the orthogonality $(\ref{Orthogonality-for-Wm})$.
\end{theo}

Observe that $\widehat{\alpha}_{m+1} > 0$, $m \geq 1$.  In fact, see Theorem \ref{Thm-ChainSequence}, that  more can be said about these coefficients if the measure  $\psi$ is such that $\int_{-1}^{1}(1-x^2)^{-1/2}d\psi(x)$ exists.

\setcounter{equation}{0}
\section{Quadrature rules associated with  $\mathcal{W}_m$  }

In order to be able to obtain the quadrature rules based on the zeros of $\mathcal{W}_m$ we first present the following theorem.

\begin{theo} \label{Thm-Division}
Let $m \geq 1$ and let $\{\mathcal{W}_m\}$ be as defined before.

Given any function $\mathcal{E} \in \Omega_{4m-1}$ there exists a function $\mathcal{F} \in \Omega_{2m-1}$ and a function $\mathcal{G} \in \Omega_{2m-1}$ such that
\[
   \mathcal{E}(x) = \mathcal{F}(x)\,\mathcal{W}_{2m}(x) + \mathcal{G}(x).
\]
Likewise, given any $\mathcal{E} \in \Omega_{4m}$ there exists a function $\mathcal{F} \in \Omega_{2m-1}$ and a function $\mathcal{G} \in \Omega_{2m}$ such that
\[
   \mathcal{E}(x) = \mathcal{F}(x)\,\mathcal{W}_{2m+1}(x) + \mathcal{G}(x).
\]

\end{theo}

\noindent {\bf Proof}. We give a proof of the latter formula.
Let $P$, $Q$, $R$ and $K_{2m+1}$ be the respective self-inversive  polynomials associated with the functions $\mathcal{E}$, $\mathcal{F}$, $\mathcal{G}$ and $\mathcal{W}_{2m+1}$. Then one needs to prove that given the self inversive-polynomial $P$ of degree at most $4m$ there exists a self-inversive polynomial $Q$ of degree at most $2m-1$ and a self-inversive polynomial $R$ of degree at most $2m$ such that
\begin{equation} \label{Eq-Division}
   P(z) = Q(z)\,K_{2m+1}(z) + z^m R(z).
\end{equation}
With the self-inversive property, we can write
\[
  \begin{array}l
   \dsp P(z) = \sum_{j=0}^{4m} p_j\, z^{j} \ = \ \sum_{j=0}^{2m-1} p_j\, z^{j} \ \, + \, p_{2m}\, z^{2m} + \sum_{j=2m+1}^{4m} \overline{p}_{4m-j}\, z^{j}, \\[3ex]
   \dsp Q(z) = \sum_{j=0}^{2m-1} q_j\, z^{j} \, = \ \sum_{j=0}^{m-1} q_j\, z^{j} \  + \sum_{j=m}^{2m-1} \overline{q}_{2m-1-j}\, z^{j}, \\[3ex]
   \dsp R(z) = \sum_{j=0}^{2m} r_j\, z^{j} \ = \ \sum_{j=0}^{m-1} r_j\, z^{j} \ \, + \, r_{m}\, z^{m} + \sum_{j=m+1}^{2m} \overline{r}_{2m-j}\, z^{j}
  \end{array}
\]
and
\[
      K_{2m+1}(z) = \sum_{j=0}^{2m+1} k_j^{(2m+1)}\, z^{j} \ = \ \sum_{j=0}^{m} k_j^{(2m+1)}\, z^{j} \  + \sum_{j=m+1}^{2m+1} \overline{k_{2m+1-j}^{(2m+1)}}\, z^{j} .
\]
Note that $p_{2m}$ and $r_{m}$ are both real.

Comparing the  coefficients of the potentials $1,z, z^2, \ldots, z^{2m}$ on both sides of (\ref{Eq-Division}) we have

\[
 \begin{array}l
  q_0 = p_0/k_0^{(2m+1)} \\[2ex]
  \dsp q_j = \big[p_j - \sum_{l=0}^{j-1} q_{l}\,k_{j-l}^{(2m+1)}\big]/k_{0}^{(2m+1)}, \quad j=1, 2, \ldots, m-1, \\[2ex]
  \dsp r_{j} = p_{m+j} - \sum_{l=0}^{m+j} q_{l}\,k_{m+j-l}^{(2m+1)}, \quad j=0, 1, \ldots, m-1, \\[2ex]
 \end{array}
\]
and
\[
    r_{m} = p_{2m} - \sum_{l=0}^{2m-1} q_{l}\,k_{2m-l}^{(2m+1)} = p_{2m} - \sum_{l=0}^{m-1}[q_{l}\,k_{2m-l}^{(2m+1)}+\overline{q_{l}\,k_{2m-l}^{(2m+1)}}].
\]
From Theorem \ref{Thm-Leading-Coeff-SIP} since  $k_{0}^{(2m+1)} \neq 0$ the above formulas are well defined. Thus, with the further observation that equalities in the higher potentials lead to the same results,  the existence of $Q$ and $R$ in (\ref{Eq-Division}) is verified.

The proof of the first formula of the theorem is similar. \hfill \qed

We now look at the interpolatory type quadrature rule at the zeros of $\mathcal{W}_m$. For the notion of polynomial interpolatory quadrature rules we cite  \cite{Gautschi-book-2004}.

First we denote the zeros of $\mathcal{W}_m$ by $x_k^{(m)}$, $k=1,2, \ldots, m$ and let $z_k^{(m)} = e^{i2\arccos(x_k^{(m)})}$, $k=1,2, \ldots, m$.  Hence, if $\mathcal{G} \in \Omega_{m-1}$ then  by Theorem \ref{Thm-Interpolation} we can write
\[
      \mathcal{G}(x) = \sum_{k=1}^{m} \mathcal{L}_{m,k}(x)\,\mathcal{G}(x_k^{(m)}),
\]
where $\mathcal{L}_{m,k} \in \Omega_{m-1}$ are functions associated with $\mathcal{W}_m$, give by
\[
   \mathcal{L}_{m,k}(x) = z^{-(m-1)/2} \ (z_k^{(m)})^{(m-1)/2}\!\!
   \prod_{\begin{array}c
           l=1 \\
           l \neq k
          \end{array}}^{m}
   \frac{z-z_l^{(m)}}{z_k^{(m)}-z_l^{(m)}}, \qquad k=1,2, \ldots, m,
\]
with $x = \cos(\theta/2)$ and $z = e^{i\theta}$. \ Consequently, if
\[
     \widehat{\lambda}_{k}^{(m)} = \int_{-1}^{1} \mathcal{L}_{m,k}(x)\, d\psi(x) \quad \mbox{and} \quad \widetilde{\lambda}_{k}^{(m)} = \int_{-1}^{1} \mathcal{L}_{m,k}(x)\,\sqrt{1-x^2}\, d\psi(x), \quad k=1,2, \ldots, m,
\]
then, for any $\mathcal{G} \in \Omega_{m-1}$, we have the following interpolatory quadrature rules
\begin{equation} \label{Eq-Interp-Quad}
    \int_{-1}^{1} \mathcal{G}(x)\, d\psi(x) = \sum_{k=1}^m \widehat{\lambda}_{k}^{(m)}\,  \mathcal{G}(x_k^{(m)}) \quad \mbox{and} \quad \int_{-1}^{1} \mathcal{G}(x)\, \sqrt{1-x^2}\,  d\psi(x) = \sum_{k=1}^m \widetilde{\lambda}_{k}^{(m)}\,  \mathcal{G}(x_k^{(m)}).
\end{equation}
Clearly, these quadrature rules hold for any distinct set of points $x_k^{(m)}$, $k=1,2, \ldots, m$.  However, since  $x_k^{(m)}$ are the zeros of the functions $\mathcal{W}_m$ we can say more.

\begin{theo} \label{Thm-Quadrature}
Let $x_k^{(m)}$, $k=1,2, \ldots, m$ be the zeros of $\mathcal{W}_m$ and let
\[
   \lambda_{k}^{(m)} = \int_{-1}^{1} \left[\mathcal{L}_{m,k}(x)\right]^2\,\sqrt{1-x^2}\, d\psi(x), \quad k=1,2,\ldots,m.
\]
If $\mathcal{E} \in \Omega_{4m-1}$, $m \geq 1$, then
\[
    \int_{-1}^{1} \mathcal{E}(x)\, d \psi(x) = \sum_{k=1}^{2m} \frac{1}{\sqrt{1-(x_{k}^{(2m)})^2}}\,\lambda_{k}^{(2m)}\,  \mathcal{E}(x_k^{(2m)})
\]
and if $\mathcal{E} \in \Omega_{4m}$, $m \geq 0$,  then
\[
    \int_{-1}^{1} \mathcal{E}(x)\,  \sqrt{1-x^2}\, d \psi(x) = \sum_{k=1}^{2m+1} \lambda_{k}^{(2m+1)}\,  \mathcal{E}(x_k^{(2m+1)}).
\]

\end{theo}

\noindent {\bf Proof}.  To obtain the quadrature rule associated with $\mathcal{E} \in \Omega_{4m-1}$, we have from Theorem \ref{Thm-Division} that there exist $\mathcal{F} \in \Omega_{2m-1}$ and  $\mathcal{G} \in \Omega_{2m-1}$ such that $\mathcal{E}(x) = \mathcal{F}(x)\,\mathcal{W}_{2m}(x) + \mathcal{G}(x)$.  Hence,  from the orthogonality given by Corollary \ref{Coro-Orthogonality-Wm} that
\[
    \int_{-1}^{1} \mathcal{E}(x)\, d \psi(x)  =  \int_{-1}^{1} \mathcal{G}(x)\, d \psi(x).
\]
Therefore, from $\mathcal{E}(x_k^{(2m)}) = \mathcal{G}(x_k^{(2m)})$, $k=1, 2, \ldots, 2m$ and from (\ref{Eq-Interp-Quad})
\[
    \int_{-1}^{1} \mathcal{E}(x)\, d\psi(x) = \sum_{k=1}^{2m} \widehat{\lambda}_{k}^{(2m)}\,  \mathcal{E}(x_k^{(2m)}),
\]
which holds for $\mathcal{E} \in \Omega_{4m-1}$. With the choice $\mathcal{E}(x) = \sqrt{1-x^2}\,[\mathcal{L}_{2m,j}(x)]^2 \in \Omega_{2m-1}$ we then obtain that
\[
  \widehat{\lambda}_{j}^{(2m)} =  \frac{1}{\sqrt{1-(x_{j}^{(2m)})^2}}\,\lambda_{j}^{(2m)}, \quad j=1,2, \ldots, 2m.
\]

Now to obtain the quadrature rule associated with $\mathcal{E} \in \Omega_{4m}$, it follows from Theorem \ref{Thm-Division} that there exist $\mathcal{F} \in \Omega_{2m-1}$ and  $\mathcal{G} \in \Omega_{2m}$ such that $\mathcal{E}(x) = \mathcal{F}(x)\,\mathcal{W}_{2m}(x) + \mathcal{G}(x)$. This leads to the interpolatory quadrature rule
\[
    \int_{-1}^{1} \mathcal{E}(x)\, \sqrt{1-x^2}\, d\psi(x) = \sum_{k=1}^{2m+1} \widetilde{\lambda}_{k}^{(2m+1)}\,  \mathcal{E}(x_k^{(2m+1)}),
\]
which holds for $\mathcal{E} \in \Omega_{4m}$.  With the choice $\mathcal{E}(x) = \,[\mathcal{L}_{2m+1,j}(x)]^2 \in \Omega_{2m}$ we then obtain that
\[
  \widetilde{\lambda}_{j}^{(2m+1)} =  \lambda_{j}^{(2m+1)}, \quad j=1,2, \ldots, 2m+1.
\]
This completes the proof of the theorem. \hfill \qed

\setcounter{equation}{0}
\section{Connection with orthogonal polynomials on the unit circle } \label{Sec-OPUC-Connection}

From now on let  $\psi$ be a positive measure on $[-1,1]$ such that $\int_{-1}^{1}(1-x^2)^{-1/2}d\psi(x)$ exists.  Let $\mu$ be a positive measure on the unit circle that satisfy
\begin{equation} \label{Eq-Measure-Connection}
        -\sin(\theta/2)\, d \mu(e^{i\theta}) =  d \psi(x),
\end{equation}
where $x = (z^{1/2}+z^{-1/2})/2 = \cos(\theta/2)$.

Observe that the measure $\mu$ is not unique. If $\widetilde{\mu}$ is a positive measure on the unit circle such that
\[
    \int_{\mathcal{C}} f(z)\, d\widetilde{\mu}(z)  = \int_{\mathcal{C}} f(z)\,d\mu(z)  + \delta\, f(1),
\]
where $\delta$ is some nonzero constant, then (\ref{Eq-Measure-Connection}) also holds for $\widetilde{\mu}$.

Now consider the sequence of self-inversive polynomials $\{\widehat{K}_m(z)\}$ given by
\[
    e^{-im\theta/2} \widehat{K}_m(e^{i\theta}) = 2^m \widehat{\mathcal{W}}_m(x), \quad m \geq 0,
\]
where $\{\widehat{\mathcal{W}}_{m}\}$ are the normalized orthogonal functions given in Theorem \ref{Thm-Orthogonality-Normalized-Wm}.  Hence, one easily obtains from Corollary \ref{Coro3-Orthogonality-Wm} and Theorem  \ref{Thm-Orthogonality-Normalized-Wm} the following.

\begin{theo} \label{Thm-Orthogonality-Normalized-Km}
The elements of sequence of polynomials $\{\widehat{K}_m(z)\}$  satisfy
\[
    \int_{\mathcal{C}} z^{-m+s} \widehat{K}_m(z)\, (1-z) d \mu(z) = 0, \quad s=0,1, \ldots, m-1, \quad m \geq 1.
\]
Moreover,
\begin{equation} \label{TTRR-MKm-thm}
  \begin{array}{l}
  \widehat{K}_0(z) = 1,\qquad \widehat{K}_1(z) = (1+i\widehat{\beta}_{1})z + (1-i\widehat{\beta}_{1}),  \\[2ex]
  \widehat{K}_{m+1}(z) = \big[(1+i\widehat{\beta}_{m+1})z + (1-i\widehat{\beta}_{m+1})] \widehat{K}_{m+1}(z) - 4 \widehat{\alpha}_{m+1} z\, \widehat{K}_{m+1}(z), \quad m \geq 1,
  \end{array}
\end{equation}
where \  $\widehat{\beta}_{m}$, $\widehat{\alpha}_{m+1}$, $m \geq 1$ are as in Theorem \ref{Thm-Orthogonality-Normalized-Wm}.

\end{theo}

The remaining results in this section, stated mainly without any proofs, follows from recent results obtained in \cite{CosFelRan-2012}.

The polynomials $\widehat{K}_m$, $m \geq 0$, are constant multiples of the CD kernels
\[
    \mathcal{K}_{m}(z,1) = \frac{\overline{s_{m+1}^{\ast}(1)}\,s_{m+1}^{\ast}(z) - \overline{s_{m+1}(1)}\,s_{m+1}(z)} {1 - z} = \sum_{j=0}^{m} \overline{s_j(1)}\,s_j(z), \quad m \geq 0,
\]
where $s_m$, $m \geq 0$, are the orthonormal polynomials with respect to the positive measure $\mu$.

The following result can be stated  for the coefficients $\{\widehat{\alpha}_{m+1}\}_{m=1}^{\infty}$.

\begin{theo} \label{Thm-ChainSequence}
Let  the  positive measure $\psi$  on $[-1,1]$ be such that the integral $\int_{-1}^{1}(1-x^2)^{-1/2}d\psi(x)$ exists. Then the sequence of positive
numbers $\{\widehat{\alpha}_{m+1}\}_{m=1}^{\infty}$ that appear in the three term recurrence formula given in Theorem \ref{Thm-Orthogonality-Normalized-Wm} (and also in the three term recurrence formula given in Theorem \ref{Thm-Orthogonality-Normalized-Km}) is a positive chain sequence.  Moreover, this positive chain sequence is such that its maximal parameter sequence does not coincide with its minimal parameter sequence.
\end{theo}

One consequence of $\{\widehat{\alpha}_{m+1}\}_{m=1}^{\infty}$ being a positive chain sequence is that this together with (\ref{TTRR-MKm-thm}) enables one to prove  also the interlacing of the zeros of $\mathcal{W}_m$ and $\mathcal{W}_{m+1}$ (see \cite{DimRan-2013}).

The results of the above Theorem may be true even without the assumed condition for the measure $\psi$.  However, we have not been able to verify this.

Even though the polynomials $\widehat{K}_m$ are uniquely defined in terms of the measure $\psi$, we have already observed that the measure $\mu$ that satisfy (\ref{Eq-Measure-Connection}) is not unique (varying according to the size of the jump at $z=1$). Hence, with such distinct measures there exist distinct sets of OPUC.
However, with $0 \leq t < 1$, if  one defines the probability measure $\mu^{(t)}$  such that
\begin{equation} \label{Eq-ProbabilityMeasure-Connection}
 \begin{array}l
       -\sin(\theta/2)\, d \mu^{(t)}(e^{i\theta}) =  c(t) \ d \psi(x),  \\[1ex]
     \ \mbox{$\mu^{(t)}$ has a jump $t$ at $z=1$ (i.e. $\mu^{(t)}$ has a pure point of size $t$ at $z=1$)}, \\[1ex]
 \end{array}
\end{equation}
where  $x = (z^{1/2}+z^{-1/2})/2 = \cos(\theta/2)$ and $c(t)$ is the normalizing constant so that $\int_{\mathcal{C}} d\mu^{(t)}(z) = 1$, then we can say more about the associated monic OPUC $S_n^{(t)}$ and hence, also the orthonormal polynomials $s_n^{(t)}$.

Let $\{M_m\}_{m=1}^{\infty}$ be the maximal parameter sequence of the positive chain sequence $\{\widehat{\alpha}_{m}\}_{m=2}^{\infty}$. Using the value of $t$ and the sequence $\{M_m\}_{n=1}^{\infty}$, we now consider the new positive chain sequence $\{\widetilde{\alpha}_{m}\}_{m=1}^{\infty}$ given by
\[
      \widetilde{\alpha}_{1} = (1-t)M_1, \qquad \widetilde{\alpha}_{m+1} = \widehat{\alpha}_{m+1} = (1-M_m)M_{m+1}, \quad m \geq 1.
\]
It is easily verified (see \cite{Chihara-book}) that the maximal parameter sequence of the positive chain sequence $\{\widetilde{\alpha}_{m}\}_{m=1}^{\infty}$ is precisely $\{M_m\}_{m=0}^{\infty}$, with  $M_0 = t$.  Let $\{\mathfrak{m}_m^{(t)}\}_{m=0}^{\infty}$ be the  minimal parameter sequence of $\{\widetilde{\alpha}_{m}\}_{m=1}^{\infty}$.  That is,
\[
  \mathfrak{m}_0^{(t)} = 0, \qquad \mathfrak{m}_1^{(t)} = \widetilde{\alpha}_{1} = (1-t)M_1, \qquad \mathfrak{m}_{m+1}^{(t)} =  \widetilde{\alpha}_{m+1}/(1-\mathfrak{m}_m^{(t)}), \ \ m \geq 1.
\]
Then we can state the following.

\begin{theo}  Let $S_m^{(t)}$, $m \geq 0$, be the monic OPUC with respect to the measure $\mu^{(t)}$ given by $(\ref{Eq-ProbabilityMeasure-Connection})$. Then the associated Verblunsky coefficients $a_{m-1}^{(t)} = -\overline{S_m^{(t)}(0)}$, $m \geq 1$, satisfy
\[
     a_{m-1}^{(t)} = \frac{1}{\tau_{m-1}}\frac{1-2\mathfrak{m}_{m}^{(t)} - i \widehat{\beta}_m}{1 - i \widehat{\beta}_m}  \quad m \geq 1,
\]
where \ $\tau_{0} = 1$  \ and \ $\dsp \tau_m = \frac{1-i\widehat{\beta}_m}{1+i\widehat{\beta}_m}\, \tau_{m-1}$, \ $m \geq 1$.
\end{theo}

\setcounter{equation}{0}
\section{Examples } \label{Sec-Examples}

Given any measure $\psi$ on $[-1,1]$, one can easily obtain by numerical computation the coefficients in the three term recurrence formulas in Theorem \ref{Thm-Orthogonality-Wm}, hence, also information about the required orthogonal functions $\mathcal{W}_m$. However, for two reasons we like to consider the normalized orthogonal functions $\widehat{\mathcal{W}}_m$ given by Theorem \ref{Thm-Orthogonality-Normalized-Wm}.  One of these reasons is that when the measure is symmetric then these functions turn out to be monic polynomials. The other reason is because of Theorem  \ref{Thm-ChainSequence}.  \\

\begin{figure}[t]
\begin{minipage}[b]{0.48\linewidth}
\centering
\includegraphics[width=0.90\linewidth]{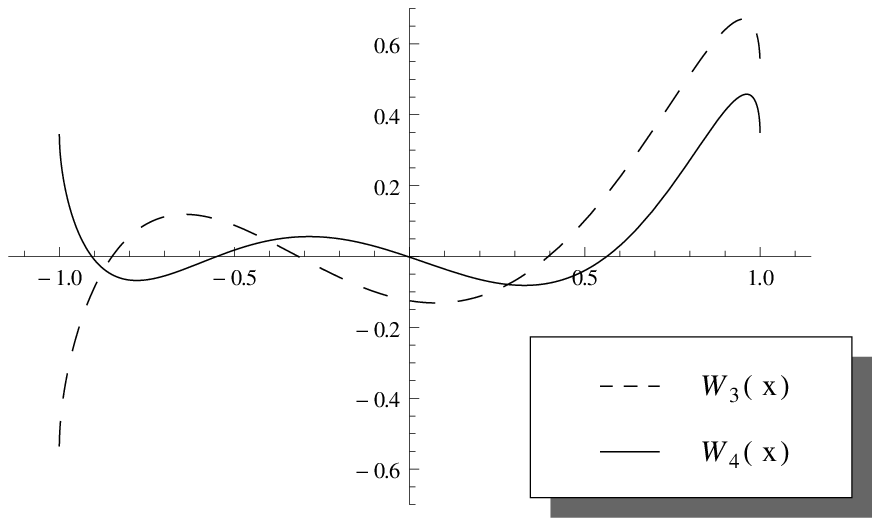}
\end{minipage} \hspace{-0.02\linewidth}
\begin{minipage}[b]{0.48\linewidth}
\centering
\includegraphics[width=0.90\linewidth]{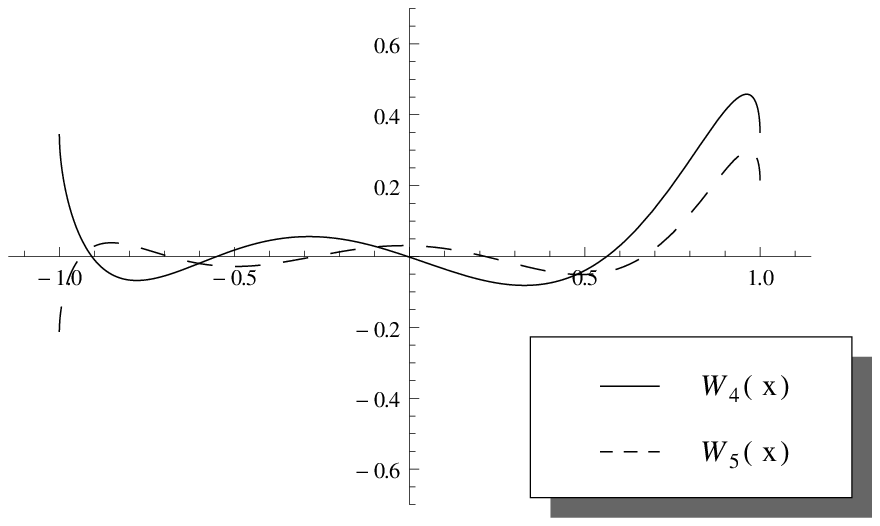}
\end{minipage}
\caption{{\small Plots of the functions $\widehat{\mathcal{W}}_{3}$, $\widehat{\mathcal{W}}_{4}$ and $\widehat{\mathcal{W}}_{5}$ when $d\psi(x) = (1-x)dx$.}}  \label{Fig1}
\end{figure}

\noindent {\bf Example 1}. \ Let $d\psi(x) = (1-x)dx$.  We obtain by numerical computation the following values for the first few $\widehat{\alpha}_m$'s and  $\widehat{\beta}_m$'s in the three term recurrence formula in Theorem \ref{Thm-Orthogonality-Normalized-Wm}. \\
\[
\begin{array}{crrrrrr}
 m & 1 \qquad & 2 \qquad & 3 \qquad & 4 \qquad & 5 \qquad & 6 \qquad \\[1ex]
 \widehat{\beta}_m & -0.4244132 & -0.3029978 & -0.2398161 & -0.2003582 & -0.1730831 & -0.1529639 \\[1ex]
 \widehat{\alpha}_m & & 0.2229581 & 0.2408213 & 0.2455306 & 0.2473987 & 0.2483152 \\[1ex]
\end{array}
\]
In the two graphs in Figure \ref{Fig1} we give, respectively,  plots of the functions  $\widehat{\mathcal{W}}_{3}$ and  $\widehat{\mathcal{W}}_{4}$ and plots of the functions $\widehat{\mathcal{W}}_{4}$ and $\widehat{\mathcal{W}}_{5}$,  separated in this way to be able see clearly the interlacing of the zeros as pointed out after Theorem \ref{Thm-ChainSequence}. Glancing at the plot of $\widehat{\mathcal{W}}_4$ it appears though this function has a zero at the origin. To be precise, the zero is very near to the origin and the value of this zero is roughly equal to  $-0.0055075$.  \\

The second example we give here is interesting from the point of view of knowing many things explicitly. \\

\noindent {\bf Example 2}. \  Let $d\psi(x) = [e^{-\arccos(x)}]^{2\eta}\, [1-x^2]^{\lambda-1} dx$, where $\eta, \lambda  \in \mathbb{R}$ and $\lambda > 1/2$. Here, we assume $\arccos(x)$  between $0$ and $\pi$.  From results given in \cite{DimIsmailRan-2012} we have
\[
     \widehat{\mathcal{W}}_{m}(x) = 2^{-m} \frac{(2\lambda)_m}{(\lambda)_m}\, e^{-im\theta/2} \, _2F_1(-m,b;\,b+\bar{b};\,1-e^{i\theta}),
\]
where $x = \cos(\theta/2)$, $b = \lambda + i \eta$ and the hypergeometric polynomial $_2F_1(-m,b;\,b+\bar{b};\,1-z)$ is self-inversive. The orthogonality of $\{\widehat{\mathcal{W}}_{m}\}$ can be explicitly written as
\[
  \begin{array}{l}
    \dsp \int_{-1}^{1} \widehat{\mathcal{W}}_{2n}(x)\,\widehat{\mathcal{W}}_{2m}(x)\, [e^{-\arccos(x)}]^{2\eta}\,[1-x^2]^{\lambda-1/2} dx = \widehat{\rho}_{2m}\,\delta_{n,m} ,  \\[2ex]
    \dsp \int_{-1}^{1} \widehat{\mathcal{W}}_{2n+1}(x)\,\widehat{\mathcal{W}}_{2m+1}(x)\, [e^{-\arccos(x)}]^{2\eta}\,[1-x^2]^{\lambda-1/2} dx =  \widehat{\rho}_{2m+1} \,\delta_{n,m},\\[2ex]
    \dsp \int_{-1}^{1} \widehat{\mathcal{W}}_{2n+1}(x)\,\widehat{\mathcal{W}}_{2m}(x)\, [e^{-\arccos(x)}]^{2\eta}\,[1-x^2]^{\lambda-1} dx = 0,
  \end{array}
\]
for $n,m = 0, 1, 2, \ldots \ $, where
\[
   \widehat{\rho}_{m} = \frac{\pi\,m! \,(\lambda+m)\, \Gamma(2\lambda+m)} {2^{2\lambda+2m-1}e^{\eta\pi}\,|\Gamma(b+m+1)|^2} \frac{1}{[(\lambda)_{m}]^2} \big[\big(\mathcal{R}e[(b)_{m}]\big)^2 + \big(\mathcal{I}m[(b)_{m}]\big)^2].
\]
Here, $\Gamma$ represents the gamma function and that  $(b)_0 = 1$ and $(b)_m = b(b+1) \cdots(b+m-1)$ for $m \geq 1$ are the Pochhammer symbols.

Moreover, in the three term recurrence formula (see Theorem \ref{Thm-Orthogonality-Normalized-Wm}) for $\{\widehat{\mathcal{W}}_{m}\}$,
\[
     \widehat{\beta}_{m} = \frac{\eta}{m+\lambda-1} \quad \mbox{and} \quad \widehat{\alpha}_{m+1} = \frac{1}{4}\frac{m (m+2\lambda-1)}{(m+\lambda-1)(m+\lambda)}, \quad m \geq 1.
\]
Observe that when $\eta = 0$, the functions $\widehat{\mathcal{W}}_{m}$ reduce to the monic ultraspherical polynomials $C_m^{(\lambda-1/2)}$.

\end{document}